\documentclass[a4paper, 12pt, oneside]{article} 
\usepackage{graphicx} 
\usepackage{amsmath} 
\usepackage{amssymb} 
\usepackage{amsfonts} 
\usepackage{amsthm} 
\usepackage{longtable} 
\usepackage{multirow}
\usepackage{color}
\usepackage{makeidx}

\makeindex
\title{Regular Polyhedra of Index Two, I
}
\author{Anthony M. Cutler\\
Northeastern University\\
Boston, MA 02115, USA\\[.05in]
{\small and} \\[.05in]
Egon Schulte\thanks{Supported by NSA-grant 
H98230-07-1-0005 and NSF-grant DMS--0856675}\\
Northeastern University\\
Boston, MA 02115, USA}
\date{ }

\newtheorem{lemma}{Lemma}[section] 
\newtheorem{theorem}[lemma]{Theorem}

\newtheorem{remark}{Remark}[section] 
    
\renewcommand{\Gamma}{\varGamma} 
\renewcommand{\epsilon}{\varepsilon}

\renewcommand{\leq}{\leqslant} 
\renewcommand{\geq}{\geqslant} 

\newcommand{\E}{\mathbb{E}^3}

\begin{document}

\maketitle

\begin{abstract}
\noindent
A polyhedron in Euclidean $3$-space $\E$ is called a {\em regular polyhedron of index $2$\/} if it is combinatorially regular but ``fails geometric regularity by a factor of $2$"; its combinatorial automorphism group is flag-transitive but its geometric symmetry group has two flag orbits. The present paper, and its successor by the first author, describe a complete classification of regular polyhedra of index $2$ in $\E$. In particular, the present paper enumerates the regular polyhedra of index $2$ with vertices on two orbits under the symmetry group. The subsequent paper will enumerate the regular polyhedra of index $2$ with vertices on one orbit under the symmetry group.
\end{abstract}

\section{Introduction}
\label{intro}

The quest for a deeper understanding of polyhedra-like structures and their symmetry in Euclidean $3$-space $\E$ has been with us since ancient times and has inspired important developments in geometry and other areas. With the passage of time, various notions of polyhedra have attracted attention and have brought to light new exciting classes of highly-symmetric figures. A key step was the radically new, graph-theoretical, or {\em skeletal\/}, approach, pioneered by Gr\"unbaum~\cite{gr1}, which introduced polyhedra as geometric (edge) graphs that are equipped with a face structure. In this vein, the notion of a regular polyhedron transcended from its origin, the five Platonic solids, to the four Kepler-Poinsot star-polyhedra, to the three Petrie-Coxeter polyhedra, finally to the more recently discovered forty-eight Gr\"unbaum-Dress polyhedra (see \cite{coxeter}, \cite{crsp}, \cite{d1}, \cite{d2}, \cite{gr1}, \cite{ordinary} and \cite[7E]{arp}). Various other classes of polyhedra with weaker symmetry properties have also been investigated and their polyhedra completely classified (see Martini~\cite{mar}, \cite{mar1} for surveys focusing on convex polyhedra). The skeletal approach has also been extended to higher-dimensional spaces, where remarkable classification results for polyhedra and polytopes have been obtained by McMullen~\cite{pm1}, \cite{pm}, \cite{pm2}, \cite{ms3}.

The present article was inspired by the study in Wills~\cite{wills} of a class of nearly regular polyhedra in $\E$. A polyhedron $P$ is called a {\em regular polyhedron of index $2$\/} if its underlying abstract polyhedron (or surface map) is combinatorially (reflexibly) regular and if its geometric symmetry group $G(P)$ has index $2$ in the combinatorial automorphism group $\Gamma(P)$ of $P$. Equivalently, $\Gamma(P)$ is flag-transitive but $G(P)$ has two orbits on the flags. Thus $P$ is combinatorially regular but ``fails geometric regularity by a factor of $2$". 

The concept of a regular polyhedron of index $2$ was introduced in Wills~\cite{wills} for orientable polyhedra with planar faces. There are exactly five such polyhedra, each with an interesting history, namely polyhedral realizations of the two pairs of dual maps $\{4,5\}_6$, $\{5,4\}_6$ of genus $4$ and $\{6,5\}_4$, $\{5,6\}_4$ of genus $9$, and of a self-dual map of type $\{6,6\}_6$ of genus $11$ (see \cite{wills} for the enumeration and Richter~\cite{rich} for figures). The polyhedra $\{5,4\}_6$ and $\{5,6\}_4$ were discovered by Hess (1878), Pitsch (1881) and Badoureau (1881), although the appeal at the time of discovery was more in their overall metrical shape than in the combinatorics of the underlying regular surface maps (abstract polyhedra) and their geometric realizations; for the history see Coxeter~\cite[\S  6.4]{coxeter}, and for figures see \cite[Fig. 4.4a]{coxeter} and Coxeter, Longuet-Higgins \& Miller~\cite[Figs. 45 and 53]{clm}. The map $\{5,4\}_6$ of genus $4$ was apparently first discovered by Gordan~\cite{gor}; under duality and Petrie-duality, it gives rise to the other three maps of genus $4$ or $9$. The polyhedra $\{4,5\}_6$ and $\{6,5\}_4$, respectively, are obtained from the polyhedra $\{5,4\}_6$ and $\{5,6\}_4$ by polarity (see Gr\"unbaum \& Shephard~\cite{gsh}); for figures see \cite[Fig. 6.4c]{coxeter}, \cite[$De_{2}f_{2}$ on Plate XI]{cdfp}, and \cite{swirp}. The polyhedron of type $\{6,6\}_6$ is an icosahedron with self-intersecting hexagonal faces, tracing back to at least Br\"uckner~\cite[Fig. 26, Tafel VIII]{brue} (see also \cite[Fig. $Ef_{1}g_{1}$ on Plate IX]{brue}).

The present paper and its successor~\cite{cut2} describe a complete classification of the finite regular polyhedra of index $2$ in $\E$ (see also \cite{cut1}), with no assumption on planarity of faces or orientability of the underlying abstract polyhedron. Here we establish that up to similarity there are precisely $22$ infinite families of regular polyhedra of index $2$ with vertices on two orbits under the full symmetry group, where two polyhedra belong to the same family if they differ only in the relative size of the spheres containing their vertex orbits. In particular, $18$ of these are related to ordinary regular polyhedra (of index $1$). All polyhedra in these $22$ families have a face-transitive symmetry group, and are orientable, but only two have planar faces (and hence occur in \cite{wills}, as $\{4,5\}_{6}$ and $\{6,5\}_{4}$). The complete enumeration of the remaining regular polyhedra of index $2$, with vertices on one orbit, is described in the subsequent paper \cite{cut2} by the first author; there are exactly $10$ such polyhedra. No finite regular polyhedron of index $2$ can be chiral; on the other hand, there are six classes of infinite regular polyhedra of index $2$ which are chiral (see \cite{chiral1}, \cite{chiral2}).

The paper is organized as follows. In Sections~\ref{mapo} and \ref{combgeo}, respectively, we cover the basics about polyhedra, maps, and combinatorial and geometric regularity. Then in Section~\ref{symgr} we investigate the symmetry groups. A key step in our approach is to determine, for all allowable vertex configurations, the possible types of faces that the polyhedron may have. To this end, in Sections~\ref{facstab} and \ref{sha}, we discuss face stabilizers and face shapes. Finally, the enumeration of polyhedra with two vertex orbits is carried out in Section~\ref{twoorbits}. 

\section{Maps and polyhedra}
\label{mapo}

Geometric realizations of abstract combinatorial complexes have attracted a lot of attention. One specific approach begins with an abstract object, in the present case an abstract polyhedron or a map on a surface, and then studies ways of realizing its combinatorics in an explicit geometric setting (for example, see 
\cite[Ch. 5]{arp} and \cite{monw5}). However, in the present paper we often take a more direct geometric approach and explicitly define the objects as polyhedra-like structures in an ambient space. Throughout, this space will be the Euclidean $3$-space $\E$. It is then understood that these objects are faithful realizations in $\E$, in the sense of \cite[Ch. 5]{arp}, of underlying abstract polyhedra; that is, they are geometric polyhedra as described below. 

An {\it abstract polyhedron}, or an {\it abstract polytope} of rank $3$, is a partially ordered set $\mathcal P$ with a strictly monotone rank function with range $\{-1,0,1,2,3\}$.  The elements of ranks $0$, $1$, or $2$ are called the {\it vertices}, {\it edges}, or {\it faces} of $\mathcal{P}$, respectively. Moreover, $\mathcal{P}$ has a smallest element (of rank $-1$) and largest element (of rank $3$), denoted by $F_{-1}$ and $F_3$, respectively. These are the {\it improper\/} elements of $\mathcal{P}$, as opposed to the other, {\em proper\/} elements of $\mathcal{P}$. Each {\it flag} (maximal chain) of $\mathcal{P}$ contains exactly five elements including $F_{-1}$ and $F_3$. Two flags are said to be {\it adjacent} if they differ in exactly one element; they are {\em $j$-adjacent\/} if this element has rank $j$. In $\mathcal{P}$, any two flags $\Phi$ and $\Psi$ can be joined by a sequence of flags $\Phi=\Phi_0,\Phi_1,...,\Phi_k=\Psi$, all containing $\Phi \cap \Psi$, such that any two successive flags $\Phi_{i-1}$ and $\Phi_i$ are adjacent; this property is known as the {\it strong flag-connectedness} of $\mathcal{P}$. Finally, $\cal P$ has the following {\it homogeneity  property}, often called the {\it diamond condition\/}:\  whenever $F\leq G$, with $F$ of rank $j-1$ and $G$ of rank $j+1$, there are exactly two elements $H$ of rank $j$ such that $F\leq H\leq G$. In practice we often ignore the least element and largest element and identify $\mathcal{P}$ with its set of proper elements, so in particular flags are $3$-element sets consisting of a vertex, an edge and a face of $P$, all mutually incident. 

All abstract polyhedra are maps on closed surfaces in the sense of Coxeter \& Moser~\cite{cm}, and all maps on closed surfaces satisfying the above homogeneity property are abstract polyhedra. Throughout we make the blanket assumption that the polyhedron or map is finite. Then the underlying surface is necessarily compact.

Informally, a geometric polyhedron will consist of a family of vertices, edges and finite polygons, all fitting together in a way characteristic for ordinary convex polyhedra or polyhedral complexes. Recall here that a (finite) {\em polygon\/}, or simply {\em $n$-gon\/}, $(v_1, v_2, \dots, v_n)$ in Euclidean $3$-space $\E$ is a figure formed by distinct points $v_1, \dots, v_n$, together with the line segments $(v_i, v_{i+1})$, for $i = 1, \dots, n-1$, and $(v_n, v_1)$. We refer to the points as {\em vertices\/} and to the line segments as {\em edges\/} of the polygon.

More precisely, a (finite) {\em geometric polyhedron} $P$ in $\E$ consists of a finite set $V_P$ of points, called {\em vertices}, a finite set $E_P$ of line segments, called {\em edges}, connecting points of $V_P$, and a finite set $F_P$ of polygons, called {\em faces}, made up of line segments of $E_P$ such that the following conditions hold. The graph defined by $V_P$ and $E_P$, called the {\em edge graph\/} of $P$, is connected. Moreover, the vertex-figure of $P$ at every vertex of $P$ is connected; by the {\em vertex-figure\/} of $P$ at a vertex $v$ we mean the graph whose vertices are the neighbors of $v$ in the edge graph of $P$ and whose edges are the line segments $(u,w)$, where $(u, v)$ and $(v, w)$ are edges of a common face of $P$. Finally, each edge of $P$ is contained in exactly two faces of $P$. For general properties of geometric polyhedra in $\E$ we refer to \cite[Ch. 7E]{arp} as well as to \cite{gr1}, \cite{grhol} and \cite{ordinary}. 

Given a face $F$ of a geometric polyhedron we sometimes abuse notation and use the term {\em boundary\/} (or {\em face boundary\/}) of $F$ to mean the underlying edge path of $F$, that is, the polygon $F$ itself.

If $P$ is a polyhedron, we let $f_{0}=f_{0}(P)$, $f_{1}=f_{1}(P)$, and $f_{2}=f_{2}(P)$ denote its number of vertices, edges, or faces, respectively. The vector $f(P)=(f_0,f_1,f_2)$ is called the 
{\em face vector\/} of $P$. 

\section{Combinatorial versus geometric regularity}
\label{combgeo}

An abstract polyhedron $\mathcal{P}$ is ({\em combinatorially\/}) {\em regular\/} if its ({\em combinatorial\/}) {\em automorphism group\/} $\Gamma(\mathcal{P})$ is transitive on the flags. The group $\Gamma(\mathcal{P})$ of a regular polyhedron $\mathcal{P}$ has a well-behaved system of \emph{distinguished generators\/} obtained as follows.
If $\mathcal{P}$ is regular, then $\Gamma(\mathcal{P})$ is generated by involutions $\rho_0,\rho_1,\rho_2$, where $\rho_j$ maps a fixed, or \emph{base\/}, flag $\Phi$ to the flag $\Phi^j$ $j$-adjacent to $\Phi$. These generators satisfy (at least) the Coxeter-type relations
\begin{equation}
\label{relone}
\rho_{0}^{2} = \rho_{1}^{2} = \rho_{2}^{2} =
(\rho_{0}\rho_{1})^{p} = (\rho_{1}\rho_{2})^{q} = (\rho_{0}\rho_{2})^{2} = \epsilon,
\end{equation}
with $\epsilon$ the identity automorphism, where $p$ and $q$ determine the ({\em Schl\"afli}) {\em type} $\{p,q\}$ of $\mathcal{P}$. Thus $\Gamma(\mathcal{P})$ is a quotient of the Coxeter group $[p,q]$ abstractly defined by the relations in (\ref{relone}). The elements
\[ \sigma_{1}:=\rho_{0}\rho_{1}, \; \sigma_{2}:=\rho_{1}\rho_{2} \]
generate the {\em rotation subgroup\/} $\Gamma^{+}(\mathcal{P})$ of $\Gamma(\mathcal{P})$, which has index at most~$2$; the index is $2$ if and only if the underlying surface is orientable. These generators satisfy the standard relations for the rotation subgroup $[p,q]^+$ of $[p,q]$, 
\[ \sigma_{1}^{p} = \sigma_{2}^{q} = (\sigma_{1}\sigma_{2})^{2} = \epsilon, \]
but in general also other independent relations. 

A {\em Petrie polygon\/} of an abstract regular polyhedron $\mathcal{P}$ (or a corresponding geometric polyhedron) is a zigzag polygon along the edges of $\mathcal{P}$ such that any two, but no three, consecutive edges belong to a common face. The {\em Petrie-dual\/} of $\mathcal{P}$ is a regular map with the same vertices and edges as $\mathcal{P}$, obtained by replacing the faces of $\mathcal{P}$ by its Petrie polygons. The Petrie-dual of $\mathcal{P}$ is a polyhedron if and only if a Petrie polygon of $\mathcal{P}$ visits any given vertex at most once (see \cite[7B3]{arp}). The iteration of the duality and Petrie operations (taking duals and Petrie-duals) forms a family of six regular maps with the same automorphism group, where possibly some maps are isomorphic. 

Recall from \cite{cm} that $\{p,q\}_r$ denotes the regular map derived from the regular tessellation $\{p,q\}$ on the $2$-sphere or the Euclidean or hyperbolic plane by identifying any two vertices $r$ steps apart along a Petrie polygon of $\{p,q\}$. Its automorphism group $[p,q]_r$ is the quotient of $[p,q] = \langle\rho_{0},\rho_{1},\rho_{2}\rangle$ obtained by factoring out the single extra relation $(\rho_{0}\rho_{1}\rho_{2})^{r}=\epsilon$. Note that $\{p,q\}_r$ is orientable if and only if $r$ is an even integer. The family of six maps derived from $\{p,q\}_r$ under duality and Petrie-duality are $\{p,q\}_{r}$, $\{q,p\}_{r}$, $\{r,q\}_{p}$, $\{q,r\}_{p}$, $\{r,p\}_{q}$ and $\{p,r\}_{q}$. If a regular polyhedron or map of type $\{p,q\}$ has Petrie polygons of length $r$, then it is a (generally proper) quotient of $\{p,q\}_r$. 

We next describe geometric regularity. The (geometric) {\em symmetry group} $G(P)$ of a geometric polyhedron $P$ consists of all isometries of $\E$ that map $P$ to itself. This group can be viewed as a subgroup of the automorphism group $\Gamma(\mathcal{P})$ of the underlying abstract polyhedron $\mathcal{P}$. A geometric polyhedron $P$ is called ({\em geometrically\/}) {\em regular} if its symmetry group $G(P)$ is transitive on the flags of $P$. If $P$ is regular, then $G(P)$ contains a set $r_0,r_1,r_2$ of generating reflections (in points, lines or planes) with the same characteristic properties as the set of generators $\rho_0,\rho_1,\rho_2$ for $\Gamma(\mathcal{P})$; in fact, $\Gamma(\mathcal{P})$ is isomorphic to $G(P)$ under the mapping taking $\rho_j$ to $r_j$ for each $j$. Thus geometric regularity of $P$ implies the combinatorial regularity of $\mathcal{P}$, or more informally, that of $P$.

The group $G(P)$ of a regular geometric polyhedron $P$ is transitive, separately, on the vertices, edges, and faces of $P$. In particular, the faces must be finite regular polygons, either planar (convex or star-) polygons or non-planar (skew) polygons (see \cite{gr1}, \cite{ordinary}). 

The group $G(P)$ generally contains both proper and improper isometries, that is, rotations, and plane reflections or rotatory reflections. By $G^{+}(P)$ we denote the subgroup of $G(P)$ of index at most $2$ consisting of the proper isometries. 

A geometric polyhedron may have many more combinatorial than geometric symmetries. If $P$ is a geometric polyhedron with underlying abstract polyhedron $\mathcal{P}$, we call the index of $G(P)$ in $\Gamma(\mathcal{P})$ the {\em index\/} of $P$ and denote it by $i(P)$. Thus 
\[ i(P)=|\Gamma(\mathcal{P}):G(P)|=|\Gamma(\mathcal{P})|/|G(P)| . \] 

A geometric polyhedron $P$ is called a {\em regular polyhedron of index $2$\/} if the underlying abstract polyhedron $\mathcal{P}$ is combinatorially regular and if $P$ itself has index $2$. Notice that combinatorial regularity of the underlying abstract polyhedron is required as part of the definition, so in particular, $P$ has a Schl\"afli symbol $\{p,q\}$. Moreover, $q\geq 3$, as $P$ is a faithful realization of $\mathcal{P}$ (but see also Remark~\ref{rem}). 

It follows from the definition that the symmetry group of a regular polyhedron of index $2$ has exactly two orbits on the flags of $P$. In fact, since $\Gamma(\mathcal{P})$ has just one flag-orbit  and $P$ cannot be geometrically regular, $P$ must have exactly two flag-orbits. (Recall that an automorphism or symmetry of $P$ is uniquely determined by its effect on a single flag.) 

We also note the following constraints on the number of transitivity classes on vertices, edges, and faces.

\begin{lemma}
\label{transvef}
If $P$ is a regular polyhedron of index $2$, then, under action of $G(P)$, $P$ has at
most two vertex orbits, at most two edge orbits, and at most two face orbits.
\end{lemma}

\begin{proof}
Since $P$ is a regular polyhedron of index $2$, it has exactly two flag orbits under $G(P)$. Select two flags of $P$ in different orbits to represent the two flag orbits. Then, if $v$ is a vertex of $P$, choose a flag that contains $v$ and map this flag to the flag representing its orbit. It follows that $v$ is in the same vertex orbit as the vertex from the representing flag. Thus there can be at most two vertex orbits of $P$ under $G(P)$. The argument for edge orbits and face orbits is similar.
\end{proof}

Next we investigate Petrie-duality in the context of regular polyhedra of index $2$.

\begin{lemma}
\label{geopetdual}
Let $P$ be a regular polyhedron of index $2$ such that any of its Petrie polygons visits any given vertex at most once. Then the Petrie-dual, $Q$, of $P$ is also a regular polyhedron of index $2$, and $G(Q) = G(P)$.
\end{lemma}

\begin{proof}
As remarked earlier, the condition on the Petrie polygons guarantees that the Petrie-dual of the underlying abstract polyhedron $\mathcal{P}$ of $P$ is again an abstract polyhedron, $\mathcal{Q}$ (see \cite[7B3]{arp}). Hence when the Petrie-duality is carried out on the geometric polyhedron $P$ itself, we arrive at a geometric polyhedron $Q$, a realization of $\mathcal{Q}$, with the same vertices and edges as $P$. By definition, the Petrie polygons of $P$ are characterized as those edge paths such that any two, but no three, consecutive edges belong to a face. This property is invariant under all combinatorial symmetries of $P$, including those arising from $G(P)$. It follows that $G(P)$ maps faces of $Q$ to faces of $Q$. Hence $G(P)$ preserves $Q$, so that $G(P)\leq G(Q) \leq \Gamma(\mathcal{Q})=\Gamma(\mathcal{P})$. In fact, since the Petrie-dual of the Petrie-dual is the original polyhedron, we must have $G(Q) = G(P)$, which is of index $2$ in $\Gamma(\mathcal{Q})$.
\end{proof}

We frequently use the following simple group-theoretical lemma, which applies in particular to $G(P)$ and its subgroup $G^{+}(P)$.

\begin{lemma}
\label{trivlem}
Let $G$ be any group, let $H$ be a subgroup of $G$ of index $2$, and let $g,h\in G$. Then $gh\in H$ if and only if either $g,h\in H$ or $g,h\not\in H$. In particular, $g^{2}\in H$, and $g\in H$ if $g$ has odd order. 
\end{lemma}

If $P$ is a regular polyhedron of index $2$, then the square of all elements in $\Gamma(\mathcal{P})$ are actually symmetries of $P$. This is a strong property. In particular, $\sigma_{1}^{2},\sigma_{2}^{2}\in G(P)$, so if $p$ or $q$ is odd, then $\sigma_{1}\in G(P)$ or $\sigma_{2}\in G(P)$, respectively; in fact, in this case we even have $\sigma_{1}\in G^{+}(P)$ or $\sigma_{2}\in G^{+}(P)$, respectively, since elements of $G(P)$ not in $G^{+}(P)$ are necessarily reflections or rotatory reflections and hence have even period. 

We usually identify a geometric polyhedron $P$ with its underlying abstract polyhedron $\mathcal{P}$. In particular, we frequently use the notion $\Gamma(P)$ in place of $\Gamma(\mathcal{P})$.

\section{The symmetry group}
\label{symgr}

Throughout this section, let $P$ be a regular polyhedron of index $2$ and type $\{p,q\}$ with group $G(P)$. Then $G(P)$, being a finite group of isometries of $\E$, must necessarily leave a point invariant, and we may take this to be the origin, $o$. Thus $G(P)$ is a finite subgroup of $O(3)$, the orthogonal group of $\E$. By $O^{+}(3)$ we denote the rotation subgroup of $O(3)$.

We begin with the following theorem that excludes all but the full Platonic symmetry groups as candidates for the symmetry groups of regular polyhedra of index $2$. 

\begin{theorem}
\label{excludegroups}
There is no regular polyhedron of index~$2$ whose symmetry group is reducible or coincides with the rotation subgroup of the symmetry group of a Platonic solid. In particular, there are no regular polyhedra of index $2$ that are planar.
\end{theorem}

\begin{proof} 
It is straightforward to eliminate the possibility of a regular polyhedron $P$ of index $2$ to lie in a plane. In fact, then $G(P)$ is necessarily a cyclic group of order $n$ or a dihedral group of order $2n$ in this plane, and $P$ has at least $n$ vertices. (The reflection in this plane is not considered to be a symmetry here since it acts trivially on $P$.) Hence, since $G(P)$ has index $2$ in $\Gamma(P)$, we must have $n\,q\leq f_{0}\,q=|G(P)|=n$ or $2n$, so necessarily $q=2$, $f_{0}=n\, (=p)$, and $G(P)$ is dihedral. Thus $P$ has just two faces and hence is not a geometric polyhedron (see also Remark~\ref{rem}). This settles the planar case.

We next exclude the reducible groups. For a finite subgroup $H$ of $O^{+}(3)$ we define $H^{*}:=H \cup (-e)H$, where $e$ denotes the identity mapping on $\E$ and $-e$ the central inversion in the origin. Moreover, if $K$ is also a finite subgroup of $O^{+}(3)$ containing $H$ as a subgroup of index $2$, we set $K]H := H \cup (-e)(K\setminus H)$. Then $-e$ is an element of the groups of type $H^*$ but not of those of type $K]H$. For the purpose of this proof let $C_n$ denote the cyclic group of order $n$ generated by a rotation by $2\pi/n$ about the $z$-axis, and let $H_n$ denote the subgroup of $O^{+}(3)$ consisting of the rotations in $C_n$ and the half-turns about the $n$ lines of symmetry of a convex regular $n$-gon in the $xy$-plane. Recall from \cite[Ch. 2]{grove} that up to conjugacy in $O(3)$ there are seven kinds of finite reducible subgroups of $O(3)$, namely $C_n$, $C_n^*$ and $C_{2n}]C_{n}$ for $n\geq 1$ and $H_n$, $H_n^*$, $H_{n}]C_{n}$ and $H_{2n}]H_{n}$ for $n\geq 2$; their orders, respectively, are $n$, $2n$ and $2n$, and $2n$, $4n$, $2n$ and $4n$.

Now let $P$ be a regular polyhedron of index $2$. Then $f_{0}=|G(P)|/q \leq |G(P)|/3$, since $q\geq 3$. We now exploit the fact that the orbits of reducible groups are usually large compared with the group orders. Suppose $G(P)$ is reducible. Inspection of the seven types of groups shows that an orbit of a point not on the $z$-axis has size at least $n$, since all groups contain $C_n$ as a subgroup, while an orbit of a point on the $z$-axis has size $1$ or $2$. On the other hand, the vertex set of $P$ is the union of at most two  such orbits, so necessarily $n\leq f_{0}\leq |G(P)|/3$ and hence $|G(P)|\geq 3n$. This only leaves the groups $H_n^*$ and $H_{2n}]H_{n}$ as possibilities. Neither group has an orbit of size $1$, so the case $f_{0}=n+1$ cannot occur here. Now observe that an orbit of these groups of size larger than $n$ must be of size at least $2n$. Hence, if $P$ has more than $n+2$ vertices, then $G(P)$ has at least one vertex orbit of size $2n$ or two vertex orbits of size $n$; in every case, $2n\leq f_{0}\leq |G(P)|/3$, which is impossible. The possibilities that $f_{0}=n+2$ can also be rejected, since then $f_{0}$ does not divide the group order, $4n$, when $n\neq 6$; the case $n=6$ can be settled directly (then $q=3$, but the vertices in the orbit of size $2$ have valency $6$). It follows that $P$ has exactly $n$ vertices, and since $C_n$ is a subgroup of the respective group, these vertices must  be those of a convex regular $n$-gon in a plane through $o$, the center of $P$. Hence $P$ lies in a plane, which is impossible. Thus $P$ cannot have a reducible symmetry group.

Finally, we refer to \cite{cut1} for a proof that the rotation subgroups of the Platonic symmetry groups can also not occur as symmetry groups of regular polyhedra $P$ of index $2$. 
\end{proof}

\begin{remark}
\label{rem}
The proof of Theorem~\ref{excludegroups} naturally leads us to investigate the degenerate case when $q=2$. Recall that a {\em dihedron\/} is a ``complex" obtained from two copies of a (planar or non-planar) polygon, the {\em faces\/} of the dihedron, by identifying them along their boundary. This is not a geometric polyhedron in our sense. Still, a dihedron is a complex of type $\{p,2\}$ for some $p\geq 3$. If the faces are regular (convex, or star, or skew) $p$-gons, then the symmetry group is isomorphic to $D_p$, the automorphism group is isomorphic to $D_{p}\times C_2$, and the dihedron is a (planar or non-planar) ``regular  complex of index $2$". Conversely, if a dihedron $P$ of type $\{p,2\}$ is a ``regular complex of index $2$", then its faces are regular polygons; in fact, every symmetry of $P$ must be a symmetry of the underlying polygon that yields the faces, and since $P$ has exactly $2p$ symmetries, the polygon has maximum possible symmetry and hence is regular. Thus the degenerate case when $q=2$ is associated with the family of dihedra that are based on regular polygons.
\end{remark} 

From now on, let $P$ be a regular polyhedron of index $2$ with an irreducible group $G(P)$, so in particular $o$ is the only invariant point of $G(P)$. Since the rotation subgroups of the Platonic symmetry groups were  excluded, this only leaves the full Platonic symmetry groups $[3,3]$, $[3,4]$ and $[3,5]$ as possibilities for $G(P)$ (see \cite[Ch. 2]{grove}). The specific nature of these groups then enables us to determine the structure of the vertex orbits under $G(P)$. The following three lemmas detail the possibilities; the first lemma is well-known. 

\begin{lemma}
\label{fundreg}
Let $S$ be a Platonic solid of type $\{r,s\}$, let $\{v, e, f\}$ be a flag of $S$, and let $T$ be the triangle whose vertices are $v$, the center of $e$, and the center of $f$. \\[-.22in]
\begin{itemize}
\item[(a)] 
Then $T$ is a fundamental region for the action of $G(S)$ on the boundary of $S$. In particular, the orbit of every point on the boundary of $S$ under $G(S)$ meets $T$ in exactly one point, and no two points of $T$ are equivalent under $G(S)$.\\[-.3in]
\item[(b)] There are seven types of orbits of boundary points of $S$ under $G(S)$, the type
depending on whether its representative point in $T$ lies at one of the three vertices
of $T$, in the relative interior of one of the edges of $T$, or in the relative interior of
$T$. Accordingly, the size of the orbit is $f_0$, $f_1$, $f_2$, or $rf_2$, $2f_1$, $rf_2$, or 
$2rf_2 = |G(S)|$, respectively, where $f_0, f_1, f_2$ are the number of vertices, edges, and faces of $S$. (Note that $rf_{2}=2f_{1}$.)
\end{itemize}
\end{lemma}

\begin{lemma}
\label{vertexsets}
If $P$ is a regular polyhedron of index $2$ with $G(P) = G(S)$, where $S$ is a Platonic solid, then each vertex orbit of $P$ under $G(P)$ is of one of the seven types described in Lemma~\ref{fundreg}, and the vertex set of $P$ is the disjoint union of one or two such orbits. 
\end{lemma}

\begin{proof}
This follows from Lemma~\ref{fundreg} if we regard $T$ as a simplicial cone emanating from the center of $S$, rather than a triangle. This cone is a fundamental region for the action of $G(S)$ on $\E$, and, up to scaling, we may choose representative vertices for each vertex orbit in one of seven ways in this cone. By Lemma~\ref{transvef}, $P$ has at most two vertex orbits.
\end{proof}

Lemma~\ref{vertexsets} enables us to determine the structure of the vertex-orbits of a regular polyhedron $P$ of index $2$. These polyhedra have full tetrahedral, octahedral, or icosahedral symmetry, and have $12$, $24$, or $60$ edges. Moreover, $qf_{0}=pf_{2}=2f_{1}$ and $3\leq q < f_0$. There are seven distinguished areas of the fundamental region $T$ that could contain a vertex. Four of these --- the three edges and the interior --- yield too many vertices to be considered. The other three produce vertices corresponding to the vertices, faces, or edges of the tetrahedron, octahedron, or icosahedron. For the vertex-orbits we thus obtain the following configurations:\\[-.3in]

\begin{tabbing}
{\it Tet\/}\={\it rahedral symmetry\/}\\
\> $q=3$, $f_{0}=8$:\quad\; \= vertices of two tetrahedra, either aligned or opposed.\\
\> $q=4$, $f_{0}=6$: \> vertices of an octahedron.\\
\\[-.15in]
{\it Octahedral symmetry\/}\\
\> $q=3$, $f_{0}=16$: \>vertices of two cubes, aligned.\\
\> $q=4$, $f_{0}=12$: \>vertices of two octahedra, aligned.\\
\> $q=4$, $f_{0}=12$: \>vertices of a cuboctahedron.\\
\> $q=6$, $f_{0}=8$:   \>vertices of a cube.\\
\\[-.15in]
{\it Icosahedral symmetry\/}\\
\> $q=3$, $f_{0}=40$: \> vertices of two dodecahedra, aligned.\\
\> $q=4$, $f_{0}=30$: \> vertices of an icosidodecahedron.\\
\> $q=5$, $f_{0}=24$: \> vertices of two icosahedra, aligned.\\
\> $q=6$, $f_{0}=20$: \> vertices of a dodecahedron.\\
\> $q=10$, $f_{0}=12$: \>  vertices of an icosahedron.
\end{tabbing}

Clearly, if $P$ has two vertex orbits we need distinguish the two cases ``aligned" and ``opposed" only if $P$ has tetrahedral symmetry.

In conclusion, we have established the following lemma.

\begin{lemma}
\label{vertplat}
The vertices of a regular polyhedron $P$ of index $2$ are located at the vertices of a single Platonic solid, a pair of aligned or opposed Platonic solids, a cuboctahedron, or an icosidodecahedron, in each case sharing the same symmetry group.
\end{lemma}

\section{Face stabilizers}
\label{facstab}

The classification of regular polyhedra of index $2$ makes continual use of the stabilizer of
a face of a polyhedron. Every such polyhedron naturally has four groups associated with it:\ the combinatorial automorphism group $\Gamma(P)$ (recall our convention to write $\Gamma(P)$ in place of $\Gamma(\mathcal{P})$), the combinatorial rotation (or even) subgroup $\Gamma^{+}(P)$ of $\Gamma(P)$ (of index at most $2$), the geometric symmetry group $G(P)$, and the geometric rotation subgroup $G^{+}(P)$ of $G(P)$ (of index at most $2$). 
Clearly,
\[ G^{+}(P)\leq G(P) \leq \Gamma(P), \]
and since $G(P)$ has index $2$ in $\Gamma(P)$, we know that $G^{+}(P)$ has index $2$ or $4$ in $\Gamma(P)$. Moreover, since $P$ is combinatorially regular, we also have $\Gamma(P) = 
\langle\rho_{0},\rho_{1},\rho_{2}\rangle$ and $\Gamma^{+}(P) = 
\langle\sigma_{1},\sigma_{2}\rangle$, where $\rho_{0},\rho_{1},\rho_{2}$ and $\sigma_{1},\sigma_{2}$ are the respective sets of distinguished generators associated with a base flag $\Phi$ of $P$ (see Section~\ref{combgeo}). 

For a face $F$ of $P$, we let $\Gamma_{F}(P)$, $\Gamma_{F}^{+}(P)$, $G_{F}(P)$, and $G_{F}^{+}(P)$ denote the stabilizer of $F$ in $\Gamma(P)$, $\Gamma^{+}(P)$, $G(P)$, or $G^{+}(P)$, respectively. Clearly, $G_{F}^{+}(P)\leq G_{F}(P) \leq \Gamma_{F}(P)$. It is well-known that, if $P$ is of type $\{p,q\}$ and $F$ is the base face in $\Phi$, then $\Gamma_{F}(P)=\langle\rho_{0},\rho_{1}\rangle\cong D_p$ and $\Gamma_{F}^{+}(P)=\langle\sigma_{1}\rangle\cong C_p$. Since $\Gamma(P)$ and $\Gamma^{+}(P)$ act transitively on the faces of $P$, any two face stabilizers are conjugate in the group and so every face stabilizer is a dihedral group $D_p$ or cyclic group $C_p$, with generating sets conjugate to  
$\rho_{0},\rho_{1}$ or $\sigma_{1}$, respectively. Note that the number of faces in the orbit of $F$ under $\Gamma(P)$, $\Gamma^{+}(P)$, $G(P)$, or $G^{+}(P)$ is just the index
\[ |\Gamma(P) : \Gamma_{F}(P)|,\; |\Gamma^{+}(P) : \Gamma_{F}^{+}(P)|,\;
|G(P):G_{F}(P)|,\; \mbox{ or } |G^{+}(P):G_{F}^{+}(P)| ,\]
respectively; for the first two groups, this number is just the total number of faces, $f_2$.

The following two lemmas describe the actions of $G(P)$ and $G^{+}(P)$ on the faces.

\begin{lemma}
\label{actgfa}
Let $P$ be a regular polyhedron of index $2$, and let $F$ be a face of $P$. Then \\[-.26in]
\begin{itemize}
\item[(a)] 
$G_{F}(P)$ is a subgroup of $\Gamma_{F}(P)\,(\cong D_p)$ of index $1$ or $2$;\\[-.3in]
\item[(b)] 
$G(P)$ has at most two orbits on the faces of $P$;\\[-.3in]
\item[(c)]
$G_{F}(P)=\Gamma_{F}(P)$ (that is, the index is $1$) if and only if there are two orbits of faces under $G(P)$.\\[-.3in]
\item[(d)]
$G_{F}(P)$ is of index $2$ in $\Gamma_{F}(P)$ if and only if $G(P)$ acts face transitively on $P$ (that is, there is just one face orbit under $G(P)$).
\end{itemize}
\end{lemma}

Thus there are two orbits of faces under $G(P)$ if and only if $|G_{F}(P)| = |\Gamma_{F}(P)| = 2p$ for some face $F$ (and hence all faces $F$), and there is just one orbit of faces under $G(P)$ if and only if $|G_{F}(P)| =|\Gamma_{F}(P)|/2 = p$ for some face $F$ (and hence all faces $F$). 

\begin{proof}
Let $i_{F} := |\Gamma_{F}(P) : G_{F}(P)|$. Then,
\begin{align}
f_{2} = |\Gamma(P) : \Gamma_{F}(P)| &=|\Gamma(P)|/|\Gamma_{F}(P)|  \notag\\
& =2|G(P)| / i_{F} |G_{F}(P)| = 2|G(P):G_{F}(P)| / i_{F} \leq 2f_{2} / i_F .\label{gf}
\end{align}
Hence $i_{F}\leq 2$, proving part (a). Then (\ref{gf}) implies $|G(P):G_{F}(P)|=i_{F}f_{2}/2 = f_{2}/2$ or $f_2$, according as $i_{F}=1$ or $2$. Hence there are at most two orbits of $G(P)$ on the faces, establishing part~(b). 

If $\Gamma_{F}(P) = G_{F}(P)$, then $f_{2} = 2|G(P) \!:\! G_{F}(P)|$, so $|G(P) \!:\! G_{F}(P)| = f_{2}/2$; that is, there are two face orbits under $G(P)$, each of size $f_{2}/2$. Conversely, if there are two orbits under $G(P)$, then $|G(P) \!:\!G_{F}(P)| = f_{2}/2$, so we must have $i_{F}= 1$, that is, $\Gamma_{F}(P) = G_{F}(P)$. This proves part (c). Finally, part (d) is just a restatement of part~(c). 
\end{proof}

Lemma~\ref{actgfa} has the following immediate consequences. If there are two orbits of faces under $G(P)$, then for any face $F$ of $P$ we have $G_{F}(P) = \Gamma_{F}(P) = \langle\rho_{0}^{F},\rho_{1}^{F}\rangle \cong D_p$, where $\rho_{0}^{F},\rho_{1}^{F}$, respectively, are suitable conjugates of $\rho_0,\rho_1$ associated with $F$ (more exactly, with a flag that contains $F$); that is, each combinatorial symmetry of $P$ that preserves $F$ can be realized by a geometric symmetry of $P$ that preserves $F$. In this case $\sigma_{1}^{F}:= \rho_{0}^{F}\rho_{1}^{F}\in G_{F}(P)$. On the other hand, if $G(P)$ acts face transitively, then there are two possibilities for the subgroup $G_{F}(P)$ of index $2$ in $\Gamma_{F}(P)$: either $G_{F}(P)\cong C_p$ ($p$ even or odd), whence $\sigma_{1}^{F}\in G_{F}(P)$; or $G_{F}(P)\cong D_{p/2}$, $p$ even, whence $\sigma_{1}^{F}\not\in G_{F}(P)$ but $(\sigma_{1}^{F})^{2}\in G_{F}(P)$, and $G_{F}(P)$ is one of the two subgroups of $\Gamma_{F}(P)$ given by
\begin{equation}
\label{twodi}
\langle \rho_{0}^{F},\sigma_{1}^{F}\rho_{0}^{F}(\sigma_{1}^{F})^{-1}\rangle \;\mbox{ or }\;
\langle \rho_{1}^{F},\sigma_{1}^{F}\rho_{1}^{F}(\sigma_{1}^{F})^{-1}\rangle .
\end{equation}

If $P$ has two face orbits under $G(P)$, then they alternate at every vertex. More exactly, we have the following lemma.

\begin{lemma}
\label{alttwo}
If $P$ is a regular polyhedron of index $2$ and $G(P)$ has two face orbits, then faces in different orbits alternate around each vertex of $P$.
\end{lemma}

\begin{proof}
If $v$ is any vertex and $\sigma_{2}^v$ is the distinguished combinatorial rotation about $v$ of period $q$ associated with (a flag that contains) $v$, then $(\sigma_{2}^v)^{2}\in G(P)$, since $G(P)$ is of index $2$ in $\Gamma(P)$. Hence alternating faces around $v$ are in the same orbit under $G(P)$. However, faces around $v$ that are adjacent cannot belong to the same orbit under $G(P)$; otherwise the faces around $v$, and hence those around any neighboring vertex of $v$, must all lie in the same orbit, so that all faces of $P$ must lie in this orbit, by the connectedness properties of $P$. Thus adjacent faces around $v$ must be in different orbits.
\end{proof}

Next we describe the action of the rotation subgroup $G^{+}(P)$ on the faces of~$P$.

\begin{lemma}
\label{gplusfa}
Let $P$ be a regular polyhedron of index $2$, and let $F$ be a face of $P$. Then \\[-.26in]
\begin{itemize}
\item[(a)] 
$G_{F}^{+}(P)$ is a subgroup of $\Gamma_{F}(P)\,(\cong D_p)$ of index $2$ or $4$;\\[-.3in]
\item[(b)] 
$G_{F}^{+}(P)$ is of index $2$ in $\Gamma_{F}(P)$ if and only if there are two orbits of faces under $G^{+}(P)$;\\[-.3in]
\item[(c)] 
$G_{F}^{+}(P)$ is of index $4$ in $\Gamma_{F}(P)$ if and only if $G^{+}(P)$ acts face transitively on $P$ (that is, there is just one face orbit under $G^{+}(P)$).
\end{itemize}
\end{lemma}

\begin{proof}
Let $j_{F} := |\Gamma_{F}(P) : G_{F}^{+}(P)|$.  Then 
\[ j_{F}= |\Gamma_{F}(P) : G_{F}^{+}(P)| = |\Gamma_{F}(P) : G_{F}(P)| \cdot |G_{F}(P) : G_{F}^{+}(P)|,\]
and both factors on the right are $1$ or $2$. Hence $j_{F}=1$, $2$, or $4$. 

We can reject the possibility that $j_{F}=1$. In fact, if $j_{F}=1$, then $\Gamma_{F}(P) = G_{F}^{+}(P)$, so $\rho_{0}^{F},\rho_{1}^{F}\in G_{F}^{+}(P)$, where again $\rho_{0}^{F},\rho_{1}^{F}$ are suitable conjugates of $\rho_0,\rho_1$ associated with (a flag that contains) $F$. Thus $\rho_{0}^{F}$ and $\rho_{1}^{F}$, being involutions, are half-turns, and so $F$ must lie in a plane through $o$. Moreover, $\sigma_{1}^{F}\in G_{F}^{+}(P)$, so $\sigma_{1}^{F}$ must be a rotation of period $p=3$, $4$, or $5$ (and $F$ is a $p$-gon). If $p=4$, then $G(P)$ is octahedral and the reflection in the plane through $F$ must be a symmetry of $P$; hence, since an edge cannot belong to more than two faces, any face adjacent to $F$ must lie in the same plane as $F$. If $p$ is odd and $F'$ is any face of $P$, then necessarily $\sigma_{1}(F')\in G_{F'}(P)$ and hence $\sigma_{1}(F')\in G_{F'}^{+}(P)$ (its order is odd). Now if $F'$ is adjacent to $F$, then also $j_{F'}=1$, since a suitable conjugate of $\rho_{0}^{F}$ by an element in $G_{F}^{+}(P)$ also belongs to $G_{F'}^{+}(P)$; hence $F'$ also lies in a plane through $o$, which must be the plane through $F$. In either case the faces adjacent to $F$ are coplanar with $F$, which is impossible. Thus the case $j_{F}=1$ cannot occur. This establishes part~(a).

Next observe that
\begin{equation}
\label{pfpf2}
f_{2} =|\Gamma(P)|/|\Gamma_{F}(P)|  
=4|G^{+}(P)| / j_{F} |G_{F}^{+}(P)| = 4|G^{+}(P):G_{F}^{+}(P)| / j_{F} \leq 4f_{2} / j_F .
\end{equation}
Hence, if $G_{F}^{+}(P)$ has index $2$ in $\Gamma_{F}(P)$, then $j_{F}=2$ and (\ref{pfpf2}) yields  
$f_{2} = 2|G^{+}(P) : G_{F}^{+}(P)|$, so the orbit of $F$ under $G^{+}(P)$ contains $f_{2}/2$ faces. Conversely, if the orbit of $F$ under $G^{+}(P)$ contains $f_{2}/2$ faces then necessarily $j_{F}=2$, that is, the index of $G_{F}^{+}(P)$ is $2$. Clearly, if the orbit of $F$ under $G^{+}(P)$ contains $f_{2}/2$ faces there are at least two face orbits, and by part~(a) they each are of size at least $f_{2}/2$ (because the case $j_{F}=1$ cannot occur). Hence there are exactly two face orbits if the orbit of $F$ under $G^{+}(P)$ contains $f_{2}/2$ faces. This settles part~(b).

Finally, part~(c) is just a reformulation of part~(b); in fact, by (\ref{pfpf2}), we have $j_{F}=4$ if and only if $f_{2} = |G^{+}(P) : G_{F}^{+}(P)|$. 
\end{proof}

According to Lemma~\ref{gplusfa}, the stabilizer groups $G_{F}^{+}(P)$ can only be of the following kind. As before, let $\Gamma_{F}(P)= \langle\rho_{0}^{F},\rho_{1}^{F}\rangle \cong D_p$, where $\rho_{0}^{F},\rho_{1}^{F}$, respectively, are suitable conjugates of $\rho_0,\rho_1$ associated with $F$. If $G_{F}^{+}(P)$ is of index $2$ in $\Gamma_{F}(P)$, then either $G_{F}^{+}(P)\cong C_p$ ($p$ even or odd), whence $\sigma_{1}^{F}\in G_{F}(P)$; or $G_{F}^{+}(P)\cong D_{p/2}$, $p$ even, whence $\sigma_{1}^{F}\not\in G_{F}(P)$ but $(\sigma_{1}^{F})^{2}\in G_{F}(P)$, and $G_{F}^{+}(P)$ is one of the groups in (\ref{twodi}). On the other hand, if $G_{F}^{+}(P)$ is of index $4$ in $\Gamma_{F}(P)$, then $p$ is even, and either $G_{F}^{+}(P) \cong C_{p/2}$ or $G_{F}^{+}(P) \cong D_{p/4}$ with $p\equiv 0 \bmod 4$. In the former case $\sigma_{1}^{F}\not\in G_{F}(P)$, but $(\sigma_{1}^{F})^{2}\in G_{F}(P)$. In the latter case $(\sigma_{1}^{F})^{2}\not\in G_{F}(P)$, but $(\sigma_{1}^{F})^{4}\in G_{F}(P)$, and $G_{F}^{+}(P)$ is one of four dihedral subgroups of $\Gamma_{F}(P)$ of order $p/2$. Although this latter case is a theoretical possibility, it does not actually occur in the enumeration.

In any case, with regards to the action of $G^{+}(P)$ on the faces of $P$ we can always be certain that $(\sigma_{1}^{F})^{4}\in G_{F}^{+}(P)$, but in principle it may occur that $(\sigma_{1}^{F})^{2}\not\in G_{F}^{+}(P)$.

\section{Face shapes}
\label{sha}

We suppose, as before, that $P$ is a regular polyhedron of index $2$ with vertices coincident with the vertices of a single Platonic solid, a pair of aligned or opposed Platonic solids, a cuboctahedron, or an icosidodecahedron. In this section we introduce for a face $F$ of $P$ a notation, called the face shape of $F$, that describes the way in which edges change their direction as one traverses the boundary of $F$. 

We noted earlier that if the vertices of $P$ are those of a centrally symmetric Platonic solid then no edge of $P$ can join two opposite vertices, for if so then that edge would subtend at least three faces under the rotational symmetries of the Platonic solid. The same is more generally true, as shown in the next lemma.

\begin{lemma}
\label{noanti}
No edge of any regular polyhedron $P$ of index $2$ joins vertices that are collinear with the center, $o$, of $P$.
\end{lemma}

\begin{proof}
Suppose that $P$ does have an edge, $e_0$, joining two vertices collinear with $o$, and let $F$ and $F'$ be the two faces of $P$ sharing $e_0$ as an edge. Then $e_0$ must lie on a $2$-fold rotation axis of $P$, and the corresponding half-turn must map $F$ to $F'$. Since $P$ has reflexive symmetry by Theorem~\ref{excludegroups}, $F$ and $F'$ must lie in a reflection plane, $H$, of $P$ that contains $e_0$, for otherwise $e_0$ would subtend four copies of $F$ under rotation and reflection. If $P$ had two face orbits under $G(P)$, then faces in different orbits would alternate around each vertex (see Lemma~\ref{alttwo}), and so $F$ and $F'$, being adjacent faces, would lie in different orbits; however, this is not the case since they are equivalent under a half-turn. Therefore $P$ has one face orbit under $G(P)$, so each face must lie in a plane through $o$. Thus every face of $P$ adjacent to $F$ must also lie in $H$. Since P is connected, it follows, by repetition of the argument, that $P$ must be planar. This is a contradiction since $P$ is finite.
\end{proof}

The vertices of each vertex orbit of $P$ are inscribed in a sphere centered at $o$, and we may take one sphere, $\mathbb{S}$, to be our {\it reference sphere\/}. Thus all, or one half of, the vertices of $P$ are inscribed in $\mathbb{S}$ and are the vertices of a Platonic solid, a cuboctahedron, or an icosidodecahedron, denoted by $S$. 

Suppose $\{v_1,v_2\}$ is an edge of $P$, where $v_1$ and $v_2$ are vertices of $P$. It is useful to consider the projection of that edge onto $\mathbb{S}$. The projection is the (shorter) arc of the great circle, which is the intersection of $\mathbb{S}$ with the plane spanned by $v_{1}$, $v_{2}$, and $o$. Since no projected edge of $P$ can join opposite vertices of $\mathbb{S}$, all projected edges are less than half a circumference in length. Note also that the projected image of each vertex of $P$ is a vertex of $S$, except when $S$ is a tetrahedron and the vertices of the other vertex orbit of $P$ lie on a tetrahedron opposed to $S$ (in this case the symmetry group is tetrahedral). Since the projection of edges is confined to $\mathbb{S}$, the angle between two projected images of adjoining edges of a face takes a value in the range between $0$ and $2\pi$.

We use this observation to introduce a notation that describes the shape of faces of $P$, and which exploits the existence of certain rotational symmetries of $P$. Each face, $F$, of $P$ is defined by the cyclically ordered vertices in its boundary, $\{v_{1}, v_{2}, v_{3},\ldots, v_{p}\}$. Suppose that we know the edge lengths of $F$ (we will clarify this in Section~\ref{twoorbits}), the vertex orbits that the vertices $v_i$ lie on, and the position of the vertices of $P$. We will see that then $F$ can also be specified by a notation of the form $[v_{1},v_{2}; a,b,c,d]$, called the {\em shape\/} (or {\em face shape\/}) of $F$; here, $\{v_{1},v_{2}\}$ is a {\em directed\/} edge, called the {\em starting edge\/}, of the boundary of $F$, and $a$, $b$, $c$, and $d$ are symbols, such as ``{\it right\/}" or ``{\it left\/}", that specify the change of direction from one edge to the next along the boundary of $F$ (traversed in the order $\{v_{1}, v_{2}, v_{3}, \ldots, v_p\}$) when it is projected onto $\mathbb{S}$. Thus $a$ represents the change of direction from $\{v_{1}, v_{2}\}$ to $\{v_{2}, v_{3}\}$ at vertex $v_2$ (as well as that from $\{v_{5}, v_{6}\}$ to $\{v_{6}, v_{7}\}$ at vertex $v_6$, etc.), while $b$ represents the change of direction from $\{v_{2}, v_{3}\}$ to $\{v_{3},v_{4}\}$ at vertex $v_3$, etc.). We explain below why four symbols --- $a, b, c, d$ --- suffice to describe the face. 

We require that the starting edges for face shapes of faces of $P$ all lie in the same edge orbit of directed edges under $G^{+}(P)$. This stipulation is necessary to properly compare faces of $P$; an example of such a comparison is given below.

The number of different directions in which a boundary path may continue is dependent on $P$. If there is an odd number of choices then the change of direction represented by the middle path will always be designated by $f$. Except in the case that $P$ has one vertex orbit and $S$ is a dodecahedron, this change of direction will always be ``({\em straight\/}) {\em forward\/}",  meaning that the projection of the two adjoining edges at that vertex lie in the same great circle. We shall see that, for all regular polyhedra of index $2$, we will have $q \leq 6$ if the edges have equal length and $q\leq 10$ otherwise, since the polyhedron must have tetrahedral, octahedral, or icosahedral symmetry. It will follow from this that there are at most five possible directions in which an edge path can continue at any vertex, and so the symbols needed to denote these changes of direction have at most five possible values. If three values suffice for the symbols of a specific polyhedron $P$, then we use $r$, $l$, and (if needed) $f$, to represent turns to the {\em right\/}, turns to the {\em left\/}, or moves {\em forward\/}. Note that $r$ and $l$ may each represent different angles on $\mathbb{S}$ for different polyhedra $P$, but that within any specific polyhedron $P$ the angles represented do not change. If $P$ needs symbols with more than three values, we use $hr$, $sr$, $sl$, and $hl$ (and possibly also $f$), representing ``{\em hard right\/}", ``{\em soft right\/}", ``{\em soft left\/}", and ``{\em hard left\/}".

It remains to give a justification why four symbols, like $a,b,c,d$, are sufficient to describe any particular face. This rests on the observation that, for any face $F$ of $P$, the automorphism $\sigma_{1}^{F}$ satisfies $(\sigma_{1}^{F})^{2}\in G(P)$ and hence $(\sigma_{1}^{F})^{4}\in G^{+}(P)$. Thus the automorphism that moves the vertices and edges of $F$ by four steps along $F$ can be realized by a rotational symmetry of $S$. So modulo $G^{+}(S)$, four directions (the $a,b,c,d$ above) are sufficient to specify $F$, irrespective of the number of vertices of $F$ or the number of vertex orbits of $P$. If the starting edge of $F$ is understood, we say $F$ is of {\em shape\/} [a,b,c,d]. In general this is not a unique designation; $F$ may also be any or all of $[b,c,d,a]$, $[c,d,a,b]$, or $[d,a,b,c]$, depending on the number of directed edge orbits of $P$ under $G^{+}(P)$. In addition, if we traverse the boundary of $F$ in the opposite direction we see that the four
shapes above are the same as $[d',c',b',a']$, $[a',d',c',b']$, $[b',a',d',c']$, and $[c',b',a',d']$,
where $x'$ represents the direction $x$ when traversing the boundary of $F$ in the opposite
direction. For example, if $x$ is ``{\em right\/}" or ``{\em left\/}", then $x'$ is ``{\em left\/}" or ``{\em right\/}", respectively. Note that while only the above eight face shapes are possible shapes for $F$, it is not necessarily the case that each of the eight represents the face shape of F. In particular,
each of the eight shapes has a different starting edge (recall that the starting edge is a directed edge), and two face shapes can represent the same face only if their starting edges are equivalent under $G^{+}(P)$.

\begin{figure}[hbt]
\centering
\begin{center}
\setlength{\unitlength}{.83pt}
\begin{picture}(300,180)
\put(-30,15){
\begin{picture}(140,150)
\put(0,0){\line(1,0){140}}
\put(140,0){\line(-1,2){30}}
\put(0,0){\line(1,2){45}}
\put(45,90){\line(-1,2){30}}
\put(110,60){\line(1,2){45}}
\put(15,150){\line(1,0){140}}
\put(0,0){\circle*{4}}
\put(140,0){\circle*{4}}
\put(45,90){\circle*{4}}
\put(110,60){\circle*{4}}
\put(15,150){\circle*{4}}
\put(155,150){\circle*{4}}
\put(50,0){\circle*{4}}
\put(105,150){\circle*{4}}
\put(-5,-12){$v_1$}
\put(45,-12){$v_2$}
\put(135,-12){$v_3$}
\put(117,58){$v_4$}
\put(10,157){$v_7$}
\put(150,157){$v_5$}
\put(100,157){$v_6$}
\put(28,88){$v_8$}
\put(7,5){$d$}
\put(47,5){$a$}
\put(127,5){$b$}
\put(100,58){$c$}
\put(23,139){$b$}
\put(141,139){$d$}
\put(102,139){$a$}
\put(49,88){$c$}
\end{picture}}
\put(180,15){
\begin{picture}(140,150)
\put(0,0){\line(1,0){140}}
\put(140,0){\line(-1,2){45}}
\put(0,0){\line(1,2){30}}
\put(30,60){\line(-1,2){45}}
\put(95,90){\line(1,2){30}}
\put(-15,150){\line(1,0){140}}
\put(0,0){\circle*{4}} 
\put(140,0){\circle*{4}} 
\put(30,60){\circle*{4}} 
\put(95,90){\circle*{4}}
\put(-15,150){\circle*{4}}
\put(125,150){\circle*{4}}
\put(90,0){\circle*{4}} 
\put(35,150){\circle*{4}}
\put(7,5){$b$}
\put(87,5){$a$}
\put(125,5){$d$}
\put(34,57){$c$}
\put(86,87){$c$}
\put(-8,139){$d$}
\put(112,139){$b$}
\put(32,139){$a$}
\end{picture}}
\end{picture}
\caption{A hypothetical face, with $p=8$, and its reflection image, with symbols denoting change of direction at each vertex.}
\label{figone}
\end{center}
\end{figure}
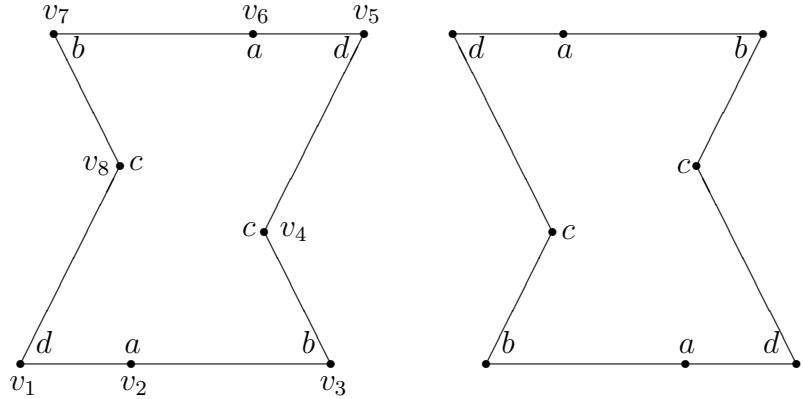

We can illustrate these observations by applying them to the hypothetical face in Figure~\ref{figone}. In this particular example, $p=8$, $a = f$ and $b = d$; nevertheless, the following statements apply generally. If we take $\{v_1,v_2\}$ as the starting edge, then the shape of the face above is $[a,b,c,d]$. However, if we have not specified the starting edge, then the shape could be any of the eight cases given in the previous paragraph. Similarly, the shape of the reflection image of the face is a cyclic permutation of either $[d,c,b,a]$ or $[a',b',c',d']$. In practice, as explained above, the (directed) starting edge is confined to a specified edge transitivity class under $G^{+}(P)$. For example, if the face in Figure~\ref{figone} has edges that can be reversed under $G^{+}(P)$, and we are using the convention that the starting edge is the shorter edge (if there are two edge lengths in $P$ and the edges of a face alternate in length), then each of the four shapes $[a,b,c,d]$, $[c,d,a,b]$, $[d',c',b',a']$, and $[b',a',d',c']$, and only those shapes, can denote the shape of the face, while its reflection image will be given by one of $[d,c,b,a]$, $[b,a,d,c]$, $[a',b',c',d']$, or $[c',d',a',b']$.

In practice it turns out that we can often do better and work with a shorter, $2$-symbol face shape notation, $[a,b]$, for a face $F$. In fact,  if we know that $(\sigma_{1}^{F})^{2}\in G^{+}(P)$, then the combinatorial automorphism that moves the vertices and edges of $F$ by two steps along $F$ can be realized by a geometric rotational symmetry of $P$, so modulo $G^{+}(P)$, two directions, $a,b$, suffice to specify $F$. In this case we use $[a,b]$ to mean $[a,b,a,b]$. We shall see in \cite{cut2} that $2$-symbol face shape notation always suffices for polyhedra of index $2$ but we have not been able to prove this a priori.

It is instructive to apply the face shape notation to some of the eighteen finite geometrically regular polyhedra (the regular polyhedra of index $1$). Since then $S$ is a Platonic solid and $G(P)=G(S)$, 
we have $(\sigma_{1}^{F})^{2} \in G^{+}(P) = G^{+}(S)$, so modulo $G^{+}(P)$ it suffices to describe a face shape by $[a,b]$. 

If $P$ is a cube, then each face of $P$ has shape $[r,r]$, where $r$ represents a turn to the right. If we traverse the boundary in the opposite direction, then the same face has shape
$[l,l]$, where $l$ represents a turn to the left. Because this is the same turn in the opposite direction, we have $l' = r$, and $r'= l$. If $P$ is the Petrie-dual of the cube, then each of its faces has shape $[r,l]$. If we traverse in the opposite direction, or choose a different starting edge, the face shape is $[l,r]$. Since the edges are equivalent modulo $G^{+}(P)$, this represents the same face shape. If $P$ is the small stellated dodecahedron $\{\frac{5}{2},5\}$, the vertices of $P$ coincide with those of an icosahedron, $S$, and $P$ has $12$ pentagram faces. Since the edge length of $P$ differs from that of $S$, the projected image on $\mathbb{S}$ (the circumsphere of $S$) of the edges of $P$ is not confined to the edges of $S$. At any vertex, $v_2$, on a face, $F$, of $P$ there are four possible directions in which to continue, relative to the edge $\{v_{1}, v_{2}\}$ of $F$, corresponding to the four vertices (other than $v_1$) which are the same distance from $v_2$ as $v_1$ is. If we call these directions in counterclockwise order $hr$, $sr$, $sl$, and $hl$ (for {\it hard right\/}, {\em soft right\/}, {\em soft left\/}, and {\it hard left\/}), we have that $[hr,hr]$ and $[hl,hl]$ represent the same shape, which is the pentagram shape of each face of $P$. For the same $S$, and edge length, the shapes $[sl,sl]$ and $[sr,sr]$ represent the faces of a great icosahedron, $\{3,\frac{5}{2}\}$.

\section{Polyhedra with two vertex orbits}
\label{twoorbits}

The regular polyhedra of index $2$ have their vertices on one or two orbits under their full symmetry group. In this section we completely enumerate the polyhedra with two vertex orbits, up to similarity. It turns out that these can be grouped into infinite families, such that the polyhedra in each family exhibit the same geometric characteristics and vary only in the relative sizes of the circumspheres of their vertex orbits. The polyhedra with one vertex orbit will be described in \cite{cut2}. 

Throughout this section, we assume that $P$ is a regular polyhedron of index $2$, of type $\{p,q\}$ (or $\{p_{P},q_{P}\}$, to emphasize $P$),  and that the vertices of $P$ lie on two orbits under $G(P)$. Then, by Lemma~\ref{vertplat} and the analysis preceding it, the vertices of $P$ are located at those of a pair of similar, aligned or opposed, Platonic solids, $S$ and $S^\diamond$, sharing with $P$ the same symmetry group and in particular the same center, $o$, and axes of rotation. Thus the vertices of one vertex orbit of $P$ are those of $S$, while the vertices of the other vertex orbit are those of $S^\diamond$. Recall our standing assumption that the vertices of $S$ lie on our reference sphere, whose size we specify shortly. Note that $S^\diamond$ is a (positively) homothetic copy of $S$ or $-S$, where the latter possibility only occurs when $S$ is the tetrahedron (and the vertices of $S$ and $-S$ give those of a cube). 

Note that the ratio $\lambda$ of the circumspheres of $S^\diamond$ and $S$ is a positive parameter naturally associated with $P$. In particular, $\lambda\neq 1$, unless $S$ is a tetrahedron and $S^\diamond=-S$ (we will see that even in this case $\lambda\neq 1$), since otherwise $S=S^\diamond$. In fact, it turns out that by varying $\lambda$ (that is, scaling one of $S$ or $S^\diamond$ while keeping the other fixed), we can produce regular polyhedra of index $2$ which have essentially the same geometric features as $P$ and are combinatorially isomorphic to $P$. Thus $P$ represents an infinite family of polyhedra, each corresponding to a different parameter value $\lambda$, and one for each $\lambda$. Polyhedra in the same family differ only in the value of $\lambda$, including when $\lambda<1$ and when $\lambda>1$. If there is no ambiguity, we use $P$ to represent both a family of polyhedra and a member of that family, relying on context to make clear which we mean.

In the present context it is convenient to rescale the reference sphere (the circumsphere of $S$) in such a way that the edges of $S$ have length $1$, and those of $S^\diamond$ have length $\lambda$, with $\lambda$ as above. 

If $S$ is a tetrahedron, then $S^\diamond$ can be inverted relative to $S$, so that it is positively homothetic to $-S$. If two polyhedra of index $2$ differ in this regard (that is, if one has $S^\diamond=\lambda S$ and the other has $S^\diamond=-\lambda S$), then we will place them into different families, even though they are combinatorially isomorphic.
\medskip

\noindent
{\bf Edge Configurations\/}\smallskip

We first establish that every edge of $P$ must join a vertex of $S$ with a vertex of $S^\diamond$.  Clearly, any edge of $P$ must either connect two vertices in $S$, or connect two vertices in $S^\diamond$, or join a vertex of $S$ with a vertex of $S^\diamond$. This property is preserved under $G(P)$, so two edges of $P$ of different kinds, and hence two flags of $P$ containing edges of different kinds, cannot be equivalent under $G(P)$. In particular, since $G(P)$ has exactly two flag orbits, edges of at most two kinds occur. Moreover, since $P$ is connected, it must contain edges joining a vertex of $S$ with a vertex of $S^\diamond$. Now consider the bipartite graph consisting of all edges of $P$ that join a vertex of $S$ to a vertex of $S^\diamond$, and their endpoints. Since $G(P)$ acts vertex-transitively on both $S$ and $S^\diamond$, its vertices are just those of $S$ and $S^\diamond$. Then either the vertices of $S$ or the vertices of $S^\diamond$ must have valency $q$ in this graph, since $P$ has edges of at most two kinds. Hence, since the graph is bipartite, all its vertices must have valency $q$. Thus every edge of $P$ joins a vertex of $S$ with a vertex of $S^\diamond$. In particular it follows that $p$, the number of vertices in a face of $P$, must be even since these vertices must alternate between the two vertex-orbits. 

For a vertex $v$ of $S$, let $v^\diamond$ denote the corresponding vertex of $S^\diamond$. So $v$ and $v^\diamond$ are collinear with $o$, the common center of $S$ and $S^\diamond$, and are on the same side of $o$ if $S$ and $S^\diamond$ are aligned. Then the observation in the previous paragraph is that every edge of $P$ is of the form $\{v_{1},v_{2}^\diamond\}$, where $v_1$ and $v_2$ are separate vertices of $S$.

Rather than using the Pythagorean length of each edge of $P$, it is more convenient to use a combinatorial metric for the length, which for simplicity we also call length. There is no ambiguity in using the same term, as the meaning defined below is always used in relation to edge length.

First suppose $P$ is such that $S$ and $S^\diamond$ are aligned (that is, $S^\diamond$ is a positively homothetic copy of $S$, namely $S^\diamond=\lambda S$). Let $\{v_{1},v_{2}^\diamond\}$ be any edge of $P$, where $v_{1},v_{2}$ are vertices of $S$ (and $v_{2}^\diamond$ is the vertex of $S^\diamond$ corresponding to $v_2$). The length of the edge $\{v_{1},v_{2}^\diamond\}$ is taken to be the combinatorial length of the shortest path from $v_1$ to $v_2$ along edges of $S$. The length of any edge of $P$ is thus an integer. 

If $P$ is such that $S$ and $S^\diamond$ are opposed (that is, $S^\diamond=-\lambda S$), then $S$ is a tetrahedron and the vertices of $S$ and $-S$ determine alternating sets of vertices of a cube. The length of any edge $\{v_{1},v_{2}^\diamond\}$ of $P$ is then taken to be the combinatorial length of the shortest path from $v_1$ to $v_{2}^\diamond$ along edges of this cube. Note that in this case the length of any edge of $P$ is $1$, since each edge of the cube joins a vertex of $S$ with a vertex of $S^\diamond$ and no edge of $P$ joins polar opposite vertices.

Next we describe the possible edge lengths for polyhedra with two vertex orbits. Note that the analogous  statement is not true for polyhedra with one vertex orbit (see \cite{cut2}).

\begin{lemma}
\label{edgtr}
Let $P$ be a regular polyhedron of index $2$ with vertices on two orbits, given by the vertex sets of $S$ and $S^\diamond$ as above. Then all edges of $P$ have the same length, which we call the {\em edge length of $P$\/}. In fact, $G(P)$ is edge-transitive. Moreover, $q_{P}=q_{S}= 3$, $4$, or $5$, where $q_P$ and $q_S$ denote the valency of a vertex of $P$ or $S$, respectively. 
\end{lemma}

\begin{proof}
The polyhedron $P$ has twice as many vertices and twice as many edges as $S$, so $P$ and $S$ have the
same vertex valency, $q$. Hence, necessarily, $q=q_{P}=q_{S}=3$, $4$, or $5$. Every $q$-fold rotational symmetry about a vertex of $S$ determines $q$ edges of $P$ of the same length emanating from this vertex. Hence all edges at a vertex of $P$ have the same length. Moreover, $G(P)$ is edge-transitive.
\end{proof}

We now examine allowable edge lengths of $P$ for each of the five possibilities for $S$, and strengthen the  statement about edge-transitivity.

\begin{lemma}
\label{allowedge}
Let $P$ be a regular polyhedron of index $2$ with vertices on two orbits, given by the vertex sets of $S$ and $S^\diamond$ as above. Then $S$ cannot be a cube. If $S$ is a tetrahedron or an octahedron, then the edge length of $P$ is $1$. If $S$ is an icosahedron, then the edge length of $P$ is either $1$ or $2$. If $S$ is a dodecahedron, then the edge length of $P$ is either $1$ or $4$. In every case, $G^{+}(P)$ is simply edge-transitive.
\end{lemma}

\begin{proof}
We first eliminate the possibility that $S$ is a cube. If $S$ is a cube, color every vertex of $S$ red or blue such that each edge of $S$ joins vertices of different colors. Since edges of $P$ go between $S$ and $S^\diamond$ and are of equal length, we have that any face of $P$ contains only red vertices of $S$ or only blue vertices of $S$. It follows that the dual edge graph of $P$ (that is, the edge graph of the dual of $P$) is not connected, since adjacent faces of $P$ must necessarily have their vertices on $S$ colored the same. Thus if $S$ is a cube, $P$ cannot be a polyhedron (but can be a compound).

If $S$ is the tetrahedron (with $S^\diamond$ aligned or opposed), octahedron, or icosahedron, the possible edge lengths are as stated, since no edge of $P$ can join a vertex of $S$ to the corresponding opposite vertex of $S^\diamond$ (by Lemma~\ref{noanti}). If $S$ is a dodecahedron, we can rule out $2$ and $3$ among the four theoretically possible choices of $1$, $2$, $3$, or $4$ for the edge length, as follows. Suppose $S$ is a dodecahedron and $P$ has edge length of either $2$ or $3$. Choose any edge, $e_0$, of $P$ and let one of its endpoints be $v_0$. There are six possible vertices that can be the other endpoint of $e_0$, and the action of the stabilizer of $v_0$ in $G(P) = G(S)$ takes $e_0$ to all of them. This contradicts $q_{P} = q_{S} = 3$. Thus the edge lengths must be $1$ or $4$ if $S$ is the dodecahedron.

Finally, the edge-transitivity of $G^{+}(P)$ follows from the following observation. Only when $S$ is a dodecahedron, do ordered pairs $\{u,v\}$ of vertices of $S$ exist that cannot be taken under the action of $G^{+}(S)=G^{+}(P)$ to every other ordered pair of vertices of $S$ the same distance apart, and then the length of $\{u,v\}$ must be $2$ or $3$. However, $2$ and $3$ were just ruled out as possible edge length. Thus $G^{+}(P)$ is edge transitive, and then also simply edge transitive since the order of $G^{+}(P)$ is just the number of edges of $P$.
\end{proof}
\medskip

\noindent
{\bf Index 1 polyhedra from index 2 polyhedra\/}\smallskip

We now exploit the geometry of the ordinary finite regular polyhedra (of index $1$) to construct new finite regular polyhedra of index $2$ with two vertex orbits. Recall that the finite regular polyhedra in $\E$ consist of the Platonic solids, the Kepler-Poinsot polyhedra, and their Petrie-duals (see \cite[Section 7E]{arp}).

\begin{lemma}
\label{index12}
Each of the $18$ finite regular polyhedra (of index $1$) determines a regular polyhedron of index $2$ (more exactly, an infinite family of such polyhedra), with its vertices on two orbits, and with nonplanar faces.
\end{lemma}

\begin{figure}[c]
\begin{center}
\includegraphics[width=16cm, height=20cm]{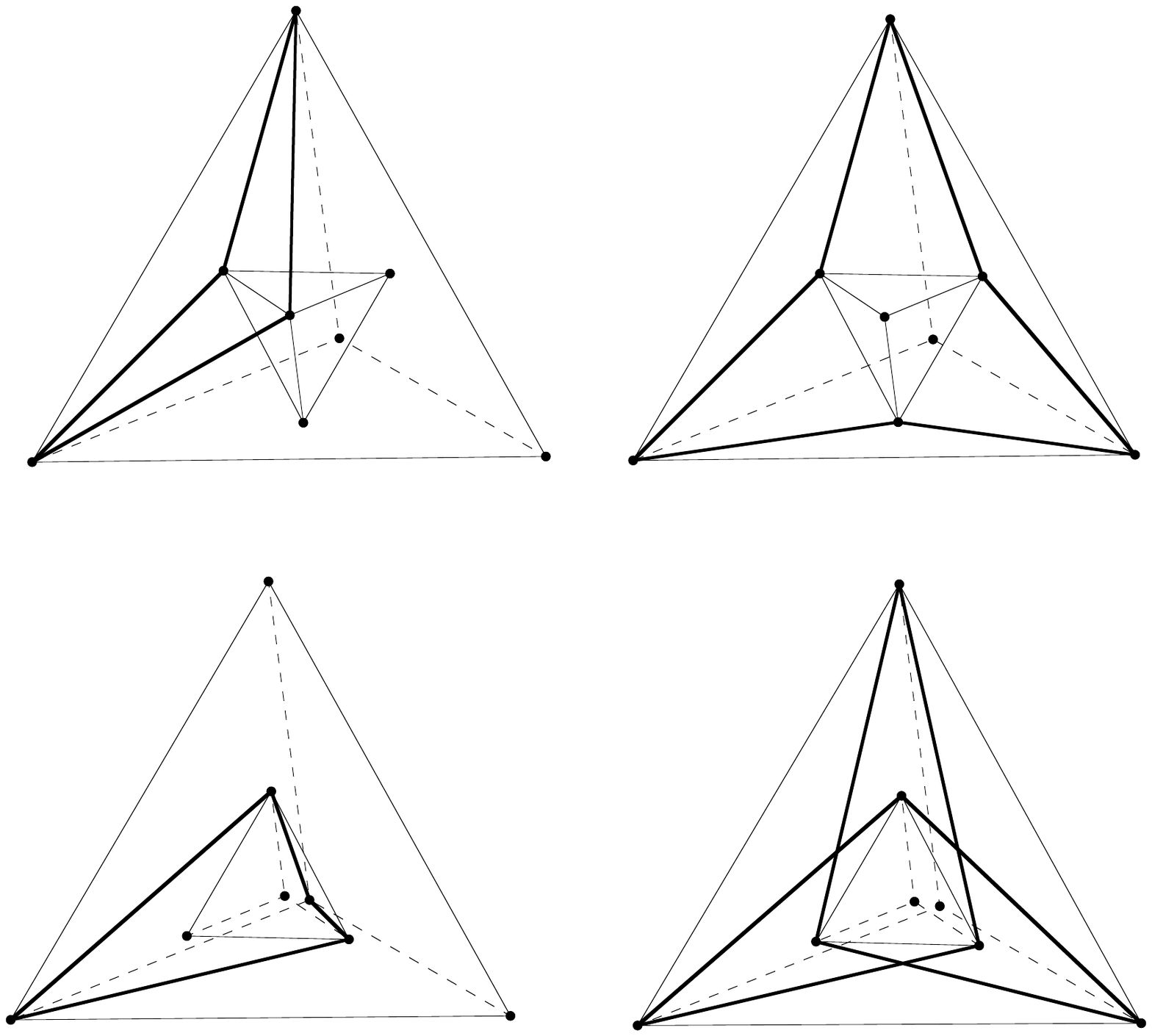}\\[-1.5in]
\caption{The four families of polyhedra with tetrahedral symmetry; from left to right, and top to bottom, derived by Lemma~\ref{index12} from the cube, the Petrie-dual of the cube, the Petrie-dual of the tetrahedron, and the tetrahedron itself. Shown is one face.}
\label{figpone}
\end{center}
\end{figure}

\begin{figure}[c]
\begin{center}
\includegraphics[width=16cm, height=19cm]{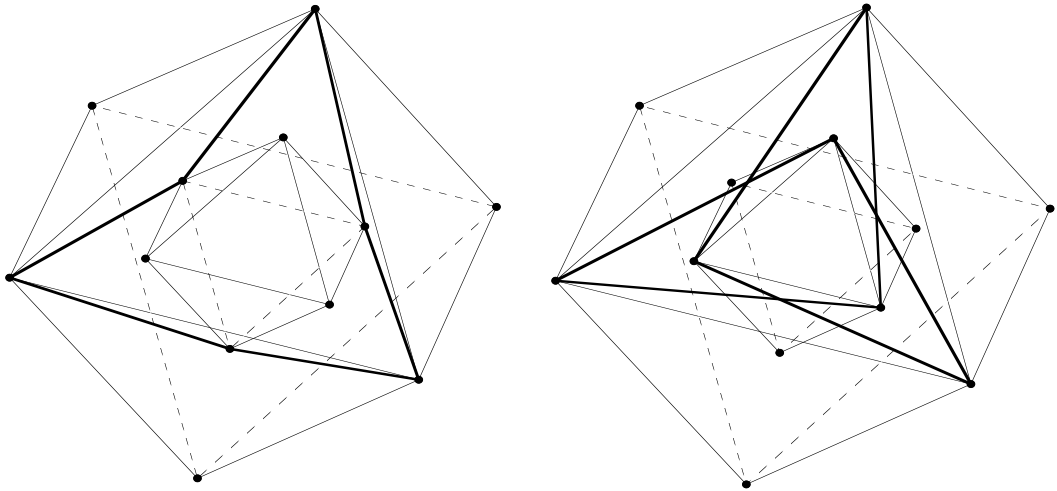}\\[-2in]
\caption{The two families of polyhedra with octahedral symmetry; from left to right, derived by Lemma~\ref{index12} from the Petrie-dual of the octahedron, and the octahedron itself. Shown is one face.}
\label{figptwo}
\end{center}
\end{figure}

\begin{figure}[c]
\begin{center}
\includegraphics[width=15cm, height=19cm]{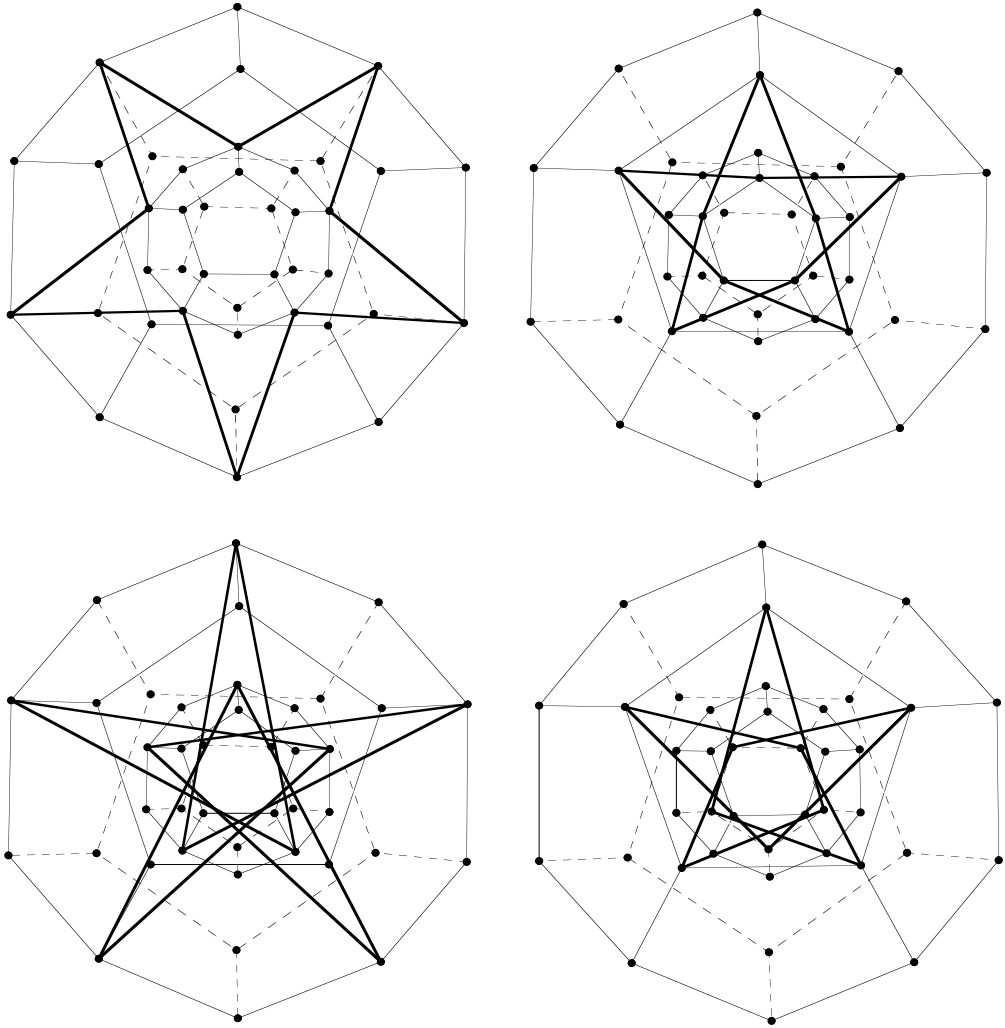}\\[-.9in]
\caption{The four families of polyhedra of type $\{10,3\}_{10}$ with icosahedral symmetry; from left to right, derived by Lemma~\ref{index12} from the Petrie-dual of the dodecahedron, the dodecahedron itself, the great stellated dodecahedron, and the Petrie-dual of the great stellated dodecahedron. Shown is one face.}
\label{figpthree}
\end{center}
\end{figure}

\begin{proof}
Let $Q$ be a finite regular polyhedron of type $\{p_{Q},q_{Q}\}$. The corresponding regular polyhedron $P$ of index $2$ and type $\{p_{P},q_{P}\}$ is constructed from $Q$ in the following manner. 

If $Q$ is the cube or its Petrie-dual, then the edge-graph of $Q$ is bipartite, and the vertices of $Q$ coincide with the vertices of two regular tetrahedra, $S$ and $S^\diamond=-S$. By rescaling $S^\diamond$ (while preserving the combinatorics of the edge graph), $Q$ is changed into a polyhedron $P$ isomorphic to $Q$ which is regular of index $2$  (see Figure~\ref{figpone}, top row). Neither of these polyhedra have planar faces: the one derived from the cube because $S$ and $S^\diamond$ have different diameters, and the one derived from the Petrie-dual of the cube because each face consists of three pairs of opposite vertices. 

If $Q$ is not the cube or its Petrie-dual, the edge-graph of $Q$ has cycles of length $3$ or $5$, and so is not bipartite. To construct $P$, take a concentric copy, $Q^\diamond$, of $Q$, positively homothetic to $Q$ but different from $Q$, and change each pair of edges between $v_0$ and $v_1$ in $Q$ and the corresponding two vertices, $v_{0}^\diamond$ and $v_{1}^\diamond$, in $Q^\diamond$ to two edges in $P$ from $v_0$ to $v_{1}^\diamond$ and from $v_{0}^\diamond$ to $v_{1}$. Both of these edges in $P$ subtend a pair of faces, each corresponding to the pair of faces of $Q$ which are subtended by $\{v_0,v_1\}$; if $q_{Q}$ is odd, then the two pairs of faces are the same.

The resulting structure is connected (the edge-graph of $Q$ is not bipartite), and is a polyhedron $P$ with $q_{P} = q_{Q}$, but with twice the number of vertices as $Q$, since all vertices of $Q$ and $Q^\diamond$ occur as vertices of $P$, and twice the number of edges as $Q$. Figure~\ref{figpone} (bottom row) and Figures~\ref{figptwo}, \ref{figpthree}, \ref{figpfive}, and~\ref{figpsix} show the vertices and a single face of each of the sixteen resulting polyhedra. If $p_{Q}$ is odd, then every face boundary of $P$ tracks the corresponding face boundary of $Q$ twice, so that $p_{P} = 2p_{Q}$; in this case, $P$ and $Q$ have the same number of faces, and every pair of faces of $P$ that are subtended by an edge of $P$ are also subtended by the opposite edge on the face boundary. If $p_{Q}$ is even, then $p_{P} = p_{Q}$, and $P$ has twice as many faces as $Q$, with each face of $Q$ corresponding to two faces of $P$.

No polyhedron in these $16$ families has planar faces. In fact, whether $p_Q$ is even or odd, any planar face of $P$ would have to lie in a plane through $o$; since the faces of $P$ are equivalent under its group, then all faces would have to lie in a plane through $o$, which is impossible. If $p_{Q}$ is odd any face of $P$ with edge $\{v_{0},v_{1}^\diamond\}$ also contains the edge $\{v_{0}^\diamond,v_{1}\}$, and these edges span a plane through $o$ that would have to contain the face if it was planar. If $p_{Q}$ is even, then $Q$ is the Petrie-dual of a Platonic solid (other than the cube) or a Kepler-Poinsot polyhedron, and so every face of $Q$ is a (skew) Petrie polygon of such a polyhedron, consisting of pairs of opposite edges; again $P$ cannot have planar faces.

It is straightforward to verify that these $18$ polyhedra are indeed regular polyhedra of index $2$ (see \cite[Appendix 1]{cut1}).
\end{proof}
\medskip

\noindent
{\bf More about face stabilizers\/}\smallskip

For polyhedra $P$ with two vertex-orbits under $G(P)$ (and hence $G^{+}(P)$), more can be said about the stabilizers $G_{F}^{+}(P)$ and $G_{F}(P)$ of a face $F$ in $G^{+}(P)$ or $G(P)$, respectively. These results complement those of Section~\ref{facstab}. In particular, the following lemma will later allow us to restrict our attention to $2$-symbol face shapes $[a,b]$. As before, suppose that the vertices of $P$ lie at those of $S$ and $S^\diamond$ and that $P$ is of type $\{p,q\}$.

\begin{lemma}
\label{stabrevisit}
If $P$ is a regular polyhedron of index $2$, with its vertices on two orbits under $G(P)$, then $G^{+}(P)$ acts transitively on the faces of $P$. Moreover, for any face $F$ of $P$, if $\sigma_{1}^{F}$ denotes the standard combinatorial rotation associated with $F$ (more exactly, with a flag containing $F$), then $(\sigma_{1}^{F})^{2}\in G^{+}(P)$ and $G_{F}^{+}(P) = \langle (\sigma_{1}^{F})^{2}\rangle \cong C_{p/2}$.
\end{lemma}

\begin{proof}
Since the edges of $P$ join vertices of $S$ and $S^\diamond$, alternate vertices of a face $F$ of $P$ necessarily lie on different vertex orbits under $G(P)$, so in particular $\sigma_{1}^{F}\not\in G(P)$ and thus  $\sigma_{1}^{F}\not\in G^{+}(P)$. Moreover, no element of $G_{F}^{+}(P)$ can interchange the endpoints of an edge of $F$ (that is, be conjugate to $\rho_{0}^{F}$ in $\Gamma_{F}(P)$), since it would necessarily have to be a half-turn about the midpoint of the edge; however, since the edge joins $S$ to $S^\diamond$, its rotation axis could not pass through $o$. Similarly, no element of $G_{F}^{+}(P)$ can fix a vertex of $F$ (that is, be conjugate to $\rho_{1}^{F}$ in $\Gamma_{F}(P)$), since it would necessarily have to be a half-turn about this vertex; this, in turn, would require $S$ to be an octahedron and the vertices of $F$ to be coplanar with $o$, which is impossible by the connectedness of $P$. It follows that $G_{F}^{+}(P)$ is a proper subgroup of $\langle\sigma_{1}^{F}\rangle$. On the other hand, 
$G_{F}^{+}(P)$ has order at least $p/2$ by Lemma~\ref{gplusfa}, so $G_{F}^{+}(P)$ is generated by $(\sigma_{1}^{F})^{2}$. Another appeal to Lemma~\ref{gplusfa} then shows that $G^{+}(P)$ acts transitively on the faces of~$P$.
\end{proof}

\begin{lemma}
\label{refface}
If $P$ is a regular polyhedron of index $2$, with its vertices on two orbits under $G(P)$, then every face $F$ of $P$ has $p/2$ planes of reflection, which each pass through opposite vertices of $F$, and are planes of reflection of $P$. In particular, $G_{F}(P) \cong D_{p/2}$, for every face $F$ of $P$.
\end{lemma}

\begin{proof}
The group $G(P)$ acts transitively on the faces of $P$, since this is already true for its subgroup $G^{+}(P)$. 
Then, by Lemma~\ref{actgfa}(d), $G_{F}(P)$ is of index $2$ in $\Gamma_{F}(P)$, for every face $F$ of $P$. As noted in the proof of Lemma~\ref{stabrevisit}, $\sigma_{1}^{F}\not\in G(P)$. Hence, the additional symmetries in the stabilizer $G_{F}(P)$ must be plane reflections, and the mirrors cannot pass through the mid-point of edges of $F$ since alternate vertices of $F$ lie on different vertex orbits. These mirrors must pass through a pair of opposite vertices of $F$, and in each case the reflection is conjugate to $\rho_{1}^{F}$. In particular, $G_{F}(P) \cong D_{p/2}$.
\end{proof}

The previous lemmas allow us to describe the faces by a single $2$-symbol face shape notation $[a,b]$, the latter meaning that every face of $P$ is of shape $[a,b,a,b]$; in fact, $G^{+}(P)$ is both face transitive and edge transitive by Lemmas~\ref{allowedge} and \ref{stabrevisit}. To make this specific, we apply it to the allowable edge configurations for $P$, set out in Lemma~\ref{allowedge}, using the convention that the (directed) starting edge goes from a vertex on $S$ to a vertex on $S^\diamond$. Since $G^{+}(P)$ is edge-transitive, it is not necessary to further specify the starting edge.

If $S$ is a tetrahedron (and the edge length is $1$) or a dodecahedron (and the edge length is $1$ or $4$), then $q = 3$ and there are only two directions in which any face boundary of $P$, when projected on the circumsphere $\mathbb{S}$ of $S$, can continue at any vertex. Call these directions $r$ and $l$ (right and left, in the standard orientation on $S$). If we traverse $F$ in the opposite direction, then $r$ and $l$ change, so that $r'= l$ and $l' = r$. The projected image of the face boundary on $\mathbb{S}$ goes along arcs corresponding to edges of $S$ if the edge length of $P$ is $1$, but if the edge length is $4$ when $S$ is a dodecahedron the projected face boundary goes along arcs corresponding to paths across faces of $S$, connecting opposite vertices of adjacent pentagonal faces.

If $S$ is an octahedron (and the edge length is $1$), $q = 4$ and we have three possible directions that a face boundary can continue at any vertex when projected on $\mathbb{S}$, namely $r$, $l$, and $f$ (where $f$ stands for ``{\em forward\/}"). Here $r'= l$, $l'= r$, and $f' = f$. 

Finally, if $S$ is an icosahedron (and the edge length is $1$ or $2$), $q = 5$ and there are four possible directions in which a face boundary may continue at any vertex when projected on $\mathbb{S}$. Calling these directions $hr$, $sr$, $sl$, and $hl$ (for {\em hard right\/}, {\em soft right\/}, {\em soft left\/}, and {\em hard left\/}), we have that $hr' = hl$, $hl'= hr$, $sr'= sl$, and $sl'= sr$. As before, the projected image of a face of $P$ onto $\mathbb{S}$ goes along arcs corresponding to edges of $S$ if the edge length of $P$ is $1$, but along arcs corresponding to paths across faces of $S$ if the edge length is $2$.
\medskip

\noindent
{\bf Enumeration\/}\smallskip

We now have sufficient restrictions on $P$ to be able to identify all regular polyhedra of index $2$ with vertices on two orbits, which we do by examining all possible face configurations for each of the allowable combinations of $S$ and edge length of $P$. Recall that Lemma~\ref{allowedge} already ruled out the possibility that $S$ is a cube.

If $S$ is a tetrahedron or a dodecahedron, $q = 3$ and so there are only two directions in which any face boundary of $P$ can continue at any vertex. So $P$ must have shape $[r,r]$, $[l,l]$, $[r,l]$, or $[l,r]$. The first two, as well as the last two, shapes represent the same face, traversed in the opposite direction. Thus there are only two possible configurations for $P$, namely $[r,r]$ and $[r,l]$, and these correspond, respectively, to the regular polyhedra of index $2$ generated as in Lemma~\ref{index12} from the following ordinary regular polyhedra:\\[-.3in]

\begin{itemize}
\item the tetrahedron, $\{3, 3\}$, and the Petrie-dual of the tetrahedron, $\{4, 3\}_3$, when $S$ is a tetrahedron and $S^\diamond$ is aligned with $S$;\\[-.3in]
\item the cube, $\{4, 3\}$, and the Petrie-dual of the cube, $\{6, 3\}_4$, when $S$ is a tetrahedron and $S^\diamond$ is opposed to 
$S$;\\[-.3in]
\item the dodecahedron, $\{5, 3\}$, and the Petrie-dual of the dodecahedron, $\{10, 3\}_5$, when $S$ is a dodecahedron and the edge length of $P$ is $1$; and\\[-.3in]
\item the great stellated dodecahedron, $\{\frac{5}{2},3\}$, and the Petrie-dual of the great stellated dodecahedron, $\{\frac{10}{3},3\}_5$, when $S$ is a dodecahedron and the edge length of $P$ is $4$.\\[-.3in]
\end{itemize}

\noindent
Note here that, if $S$ is a tetrahedron and $S^\diamond$ is opposed to $S$, the projected image of a face of $P$ on the circumsphere of $S$ takes, as its vertices, the projections of certain vertices and face centers of $S$. The face shapes $[r,r]$ and $[r,l]$ are taken relative to this system of points. This results in the cube or its Petrie-dual. It turns out that only when the vertices of the two orbits are those of opposed tetrahedra do we get the generating (index~$1$) polyhedron to have a symmetry group different from that of the resultant (index~$2$) polyhedron.

If $S$ is an octahedron, $q = 4$ and there are nine possible configurations for $P$, represented by $[r,r]$, $[r,l]$, $[r,f]$, $[l,r]$, $[l,l]$, $[l,f]$, $[f,r]$, $[f,l]$, and $[f,f]$. However, $[f,f]$ does not produce a (connected) polyhedron because all vertices in any face boundary are coplanar with $o$. Eliminating this, as well as duplicates caused by traversing faces in the reverse direction, we have that the only possible face configurations for $P$ are $[r,r]$, $[r,l]$, $[r,f]$, and $[f,r]$. Examining these, we have that $[r,r]$ produces the regular polyhedron of index $2$ generated, as in Lemma~\ref{index12}, from the octahedron, $\{3,4\}$, and $[r,l]$ produces the polyhedron generated from the Petrie-dual of the octahedron, $\{6,4\}_3$. The configuration $[r,f]$, which represents a face with $p = 4$, does not produce a polyhedron because every face boundary is of the cyclically ordered form $\{u,v,w,-v\}$, where $-v$ is the vertex opposite $v$. In fact, for every face with cyclically ordered vertices $\{u,v,w,-v\}$, there is also a face with cyclically ordered vertices $\{-u,v,w,-v\}$, with $-u$ opposite to $u$; these two faces both contain the adjoining edges $\{v,w\}$ and $\{w,-v\}$ in their face boundaries. Then these two faces completely determine the neighborhood of $w$ (topologically speaking), forcing $q=2$. Thus $[r,f]$ cannot give a polyhedron. Finally, the configuration $[f,r]$ is combinatorially equivalent to $[r,f]$ (and geometrically equivalent if we interchange $S$ and $S^\diamond$ and rescale so that the edge length of $S$ is $1$), and so does not produce a polyhedron either. Thus there are exactly two regular polyhedra of index $2$ when $S$ is an octahedron, which are the polyhedra generated, as in Lemma~\ref{index12}, from the octahedron and its Petrie-dual.

If $S$ is an icosahedron, $q = 5$ and there are now sixteen possible configurations for $P$, corresponding to the four possible directions in which a face boundary may continue at any vertex. As before, these include traversing the face boundary in the opposite direction. Eliminating these duplications leaves the ten configurations $[hr,hr]$, $[hr,hl]$, $[hl,hr]$, $[sr,sr]$, $[sr,sl]$, $[sl,sr]$, $[hr,sr]$, $[sr,hr]$, $[hr,sl]$, and $[sl,hr]$. In addition, $[hr,hl]$ and $[hl,hr]$ are geometrically equivalent, as are each of the three pairs $[sr,sl]$, $[sl,sr]$ and $[hr,sr]$, $[sr,hr]$ and $[hr,sl]$, $[sl,hr]$, all four by interchanging $S$ and $S^\diamond$ and rescaling. Thus there are six possible configurations, for families of geometrically equivalent polyhedra, namely $[hr,hr]$, $[hr,hl]$, $[sr,sr]$, $[sr,sl]$, $[hr,sr]$, and $[hr,sl]$, for both edge length of $1$ and edge length of $2$. All twelve of these configurations produce regular polyhedra of index $2$, some quite remarkable. In each case it is straightforward to check that the polyhedron is combinatorially regular. Together with the ten examples we have just described, this list completes the enumeration of all regular polyhedra of index $2$ with two vertex orbits. We further describe the twelve families of polyhedra derived from the icosahedron as follows.

Eight of these (families of) polyhedra are obtained, as in Lemma~\ref{index12}, from ordinary regular polyhedra. The $[hr,hr]$ configuration produces the polyhedron generated from the icosahedron, $\{3,5\}$, for edge length $1$, and from the small stellated dodecahedron, $\{\frac{5}{2}, 5\}$, for edge length $2$; and the $[hr,hl]$ configuration gives the polyhedra derived from the Petrie-duals of $\{3, 5\}$ and $\{\frac{5}{2}, 5\}$ for edge lengths $1$ and $2$, respectively. The $[sr,sr]$ configuration produces the polyhedron generated from the great dodecahedron, $\{5,\frac{5}{2}\}$, for edge length $1$, and from the great icosahedron, $\{3,\frac{5}{2}\}$, for edge length $2$; and the $[sr,sl]$ configuration gives the polyhedra generated from the Petrie-duals of $\{5,\frac{5}{2}\}$ and $\{3,\frac{5}{2}\}$ for edge lengths $1$ and $2$, respectively. 

\begin{figure}[c]
\begin{center}
\includegraphics[width=15cm, height=19cm]{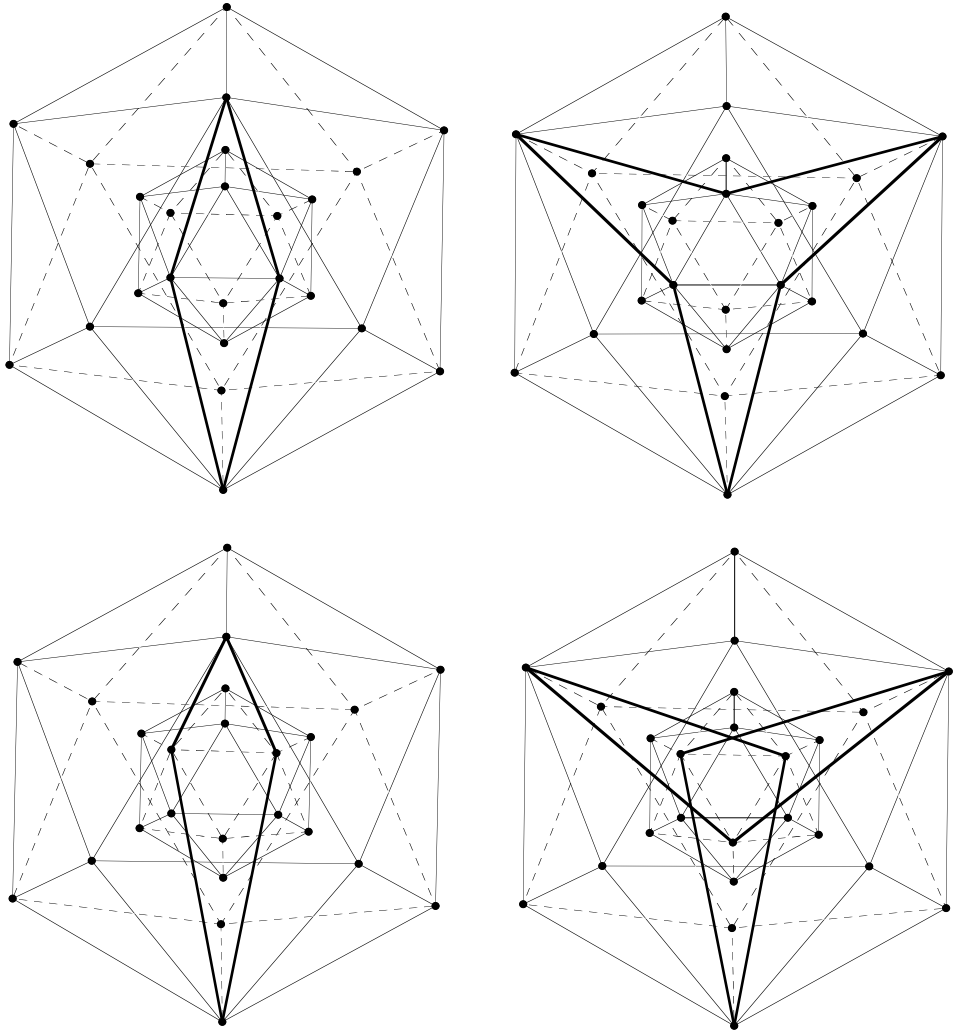}\\[-.9in]
\caption{The four exceptional families of polyhedra of types $\{4,5\}_{6}$ or $\{6,5\}_4$ with icosahedral symmetry, not derived from ordinary regular polyhedra. The two families at the top each contains one polyhedron with planar faces, obtained when $\lambda=(1+\sqrt{5})/2$ or $2+\sqrt{5}$, respectively.}
\label{figpfour}
\end{center}
\end{figure}

\begin{figure}[c]
\begin{center}
\includegraphics[width=15cm, height=19cm]{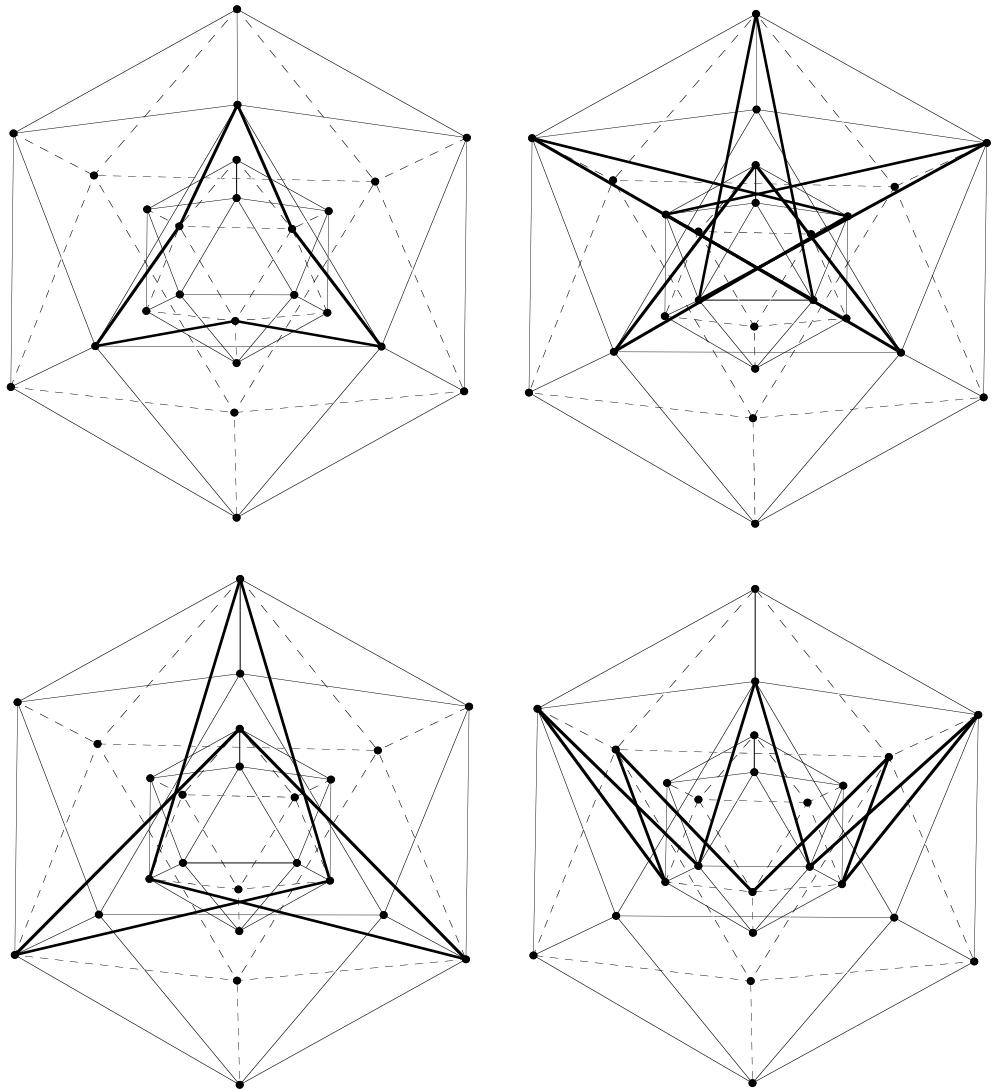}\\[-.9in]
\caption{The first set of four families of polyhedra of type $\{6,5\}_{10}$ or $\{10,5\}_6$ with icosahedral symmetry; from left to right, derived by Lemma~\ref{index12} from the Petrie-dual of the small stellated dodecahedron, the small stellated dodecahedron, the icosahedron, and the Petrie-dual of the icosahedron. Shown is one face.}
\label{figpfive}
\end{center}
\end{figure}

\begin{figure}[c]
\begin{center}
\includegraphics[width=15cm, height=19cm]{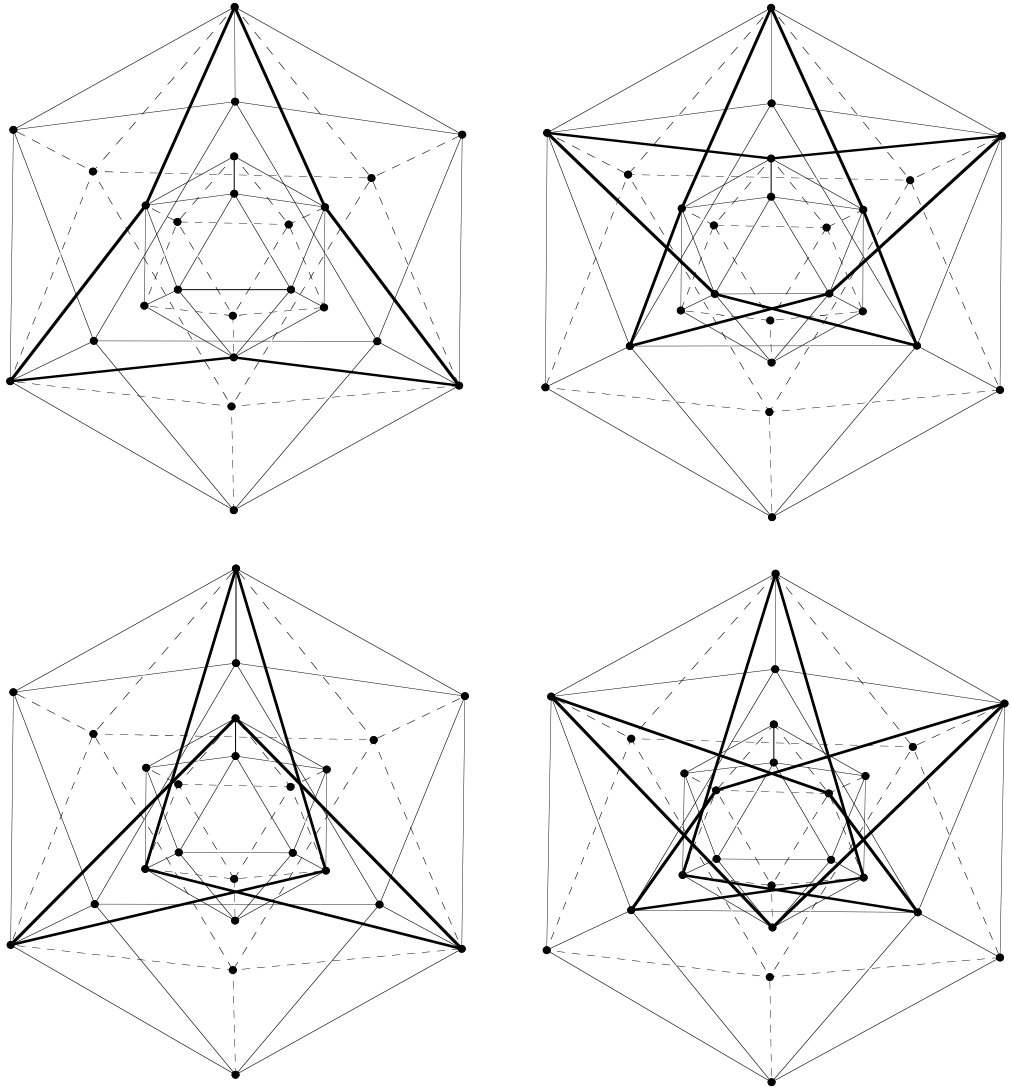}\\[-.9in]
\caption{The second set of four families of polyhedra of type $\{10,5\}_{6}$ or $\{10,5\}_6$ with icosahedral symmetry; from left to right, derived by Lemma~\ref{index12} from the Petrie-dual of the great dodecahedron, the great dodecahedron, the great icosahedron, and the Petrie-dual of the great icosahedron. Shown is one face.}
\label{figpsix}
\end{center}
\end{figure}

This list of eight is complemented by another group of four families of polyhedra (see Figure~\ref{figpfour}) which cannot be obtained by the method described in Lemma~\ref{index12}. Specifically, the $[hr,sr]$ configuration with edge length $1$ and the $[hr,sl]$ configuration with edge length $2$ determine two combinatorially equivalent orientable polyhedra of type $\{4, 5\}_6$, and the $[hr,sl]$ configuration with edge length $1$ and the $[hr,sr]$ configuration with edge length $2$ produce two combinatorially equivalent orientable polyhedra of type $\{6, 5\}_4$. These cannot be obtained from ordinary regular polyhedra by means of Lemma~\ref{index12}. In fact, if the vertices are forced into one orbit, by making $\lambda = 1$, the vertices and edges of the two orbits under $G(P)$ coincide, halving in number, but the faces do not; consequently, each edge would subtend four faces, and the resulting object would not be a polyhedron. 

We look now at whether any polyhedra in these four families have planar faces. (We already know by Lemma~\ref{index12} that each polyhedron in the other eighteen families has non-planar faces.)  For any polyhedron in either of the two families of polyhedra of type $\{4,5\}_6$, the convex hull of the vertices of any face contains $o$ if the edge length is $2$, and does not if the edge length is $1$. The vertices of any face of a polyhedron of either of the two families of polyhedra of type $\{6, 5\}_4$ lie on two parallel planes, one for each orbit, which are separated by $o$ when the edge length is $2$, but not if the edge length is $1$. Thus, polyhedra of either family with edge length $2$ do not have planar faces, but for each of the two polyhedral families with edge length $1$ there is exactly one polyhedron (denoted by a specific value of $\lambda$) which has planar faces. This occurs when $\lambda = \tau$ for the polyhedral family of type $\{4,5\}_6$, and $\lambda = 2\tau + 1$ for the polyhedral family of type $\{6,5\}_4$, where $\tau = (1+\sqrt{5})/2$ is the golden ratio. (Recall that $\lambda$ is the edge length of $S^\diamond$, and so is the ratio of the diameters of $S^\diamond$ and $S$.) Note also that ``swapping" the orbit locations $S$ and $S^\diamond$ of $P$, by changing $\lambda$ to the reciprocal of its original value, will result in a geometrically similar polyhedron if and only if $P$ has been generated, as specified in Lemma~\ref{index12}, from an ordinary regular polyhedron; polyhedra belonging to any of the four families of polyhedra just described do not maintain geometric similarity if the vertex orbits are swapped.

Thus we have established the following result.

\begin{theorem}
\label{enumthm}
There are precisely 22 infinite families of regular polyhedra of index $2$ with vertices on two orbits under the full symmetry group, where two polyhedra belong to the same family if they differ only in the relative size of the spheres containing their vertex orbits. In particular, 18 of these are related to ordinary regular polyhedra (of index 1) as described in Lemma~\ref{index12}. Of these 22 families, the polyhedra in four have tetrahedral symmetry, in two have octahedral symmetry, and in sixteen have icosahedral symmetry.  All polyhedra in these 22 families have a face-transitive symmetry group, and are orientable, but only two among them have planar faces.
\end{theorem}

\begin{table}[htb]
\centering
{\begin{tabular}{|c|c|c|c|c|c|c|}  \hline
Type & Generated by  &Face Vector  &Edge              &Face&Map &Figure \\
$\{p,q\}_r$ &Lemma~\ref{index12} from &$(f_{0},f_{1},f_{2})$&Length&Shape& of \cite{con}&\\[.05in]
\hline
\hline
$\{4,3\}_6$  &  Cube & $(8, 12, 6)$& $1$ &$[r,r]$&---&2\\
\hline
$\{6,3\}_4$  &  Petrie-dual of Cube & $(8, 12, 4)$& $1$ &$[r,l]$&---&2 \\
\hline
$\{4,3\}_6$   & Petrie-dual of  & $(8, 12, 6)$& $1$&$[r,l]$ &---&2 \\
&Tetrahedron&&&&&\\
\hline
$\{6,3\}_4$ &  Tetrahedron & $(8, 12, 4)$ & $1$ &$[r,r]$&---&2 \\
\hline
\hline
$\{6,4\}_6$  &  Petrie-dual of & $(12, 24, 8)$& $1$& $[r,l]$&$R3.4^*$ &3\\
&Octahedron&&&&&\\
\hline
$\{6,4\}_6$  &  Octahedron& $(12, 24, 8)$& $1$ &$[r,r]$ &$R3.4^*$&3 \\
\hline
\hline
$\{10,3\}_{10}$  &  Petrie-dual of & $(40,60,12)$&$1$&$[r,l]$ &$R5.2^*$&4\\
&Dodecahedron&&&&&\\
\hline
$\{10,3\}_{10}$  &  Dodecahedron& $(40,60,12)$& $1$&$[r,r]$&$R5.2^*$ &4\\
\hline
$\{10,3\}_{10}$  &  Great Stellated & $(40,60,12)$& $4$&$[r,r]$&$R5.2^*$ &4\\
&Dodecahedron&&&&&\\
\hline
$\{10,3\}_{10}$  &  Petrie-dual of Great& $(40,60,12)$ &$4$&$[r,l]$&$R5.2^*$&4 \\
&Stellated Dodecahedron&&&&&\\
\hline
$\{4,5\}_{6}$  &  ---& $(24,60,30)$ &$1$&$[hr,sr]$ &$R4.2$&5 \\
\hline
$\{6,5\}_{4}$  &  --- & $(24,60,20)$& $1$&$[hr,sl]$&$R9.16^*$ &5\\
\hline
$\{4,5\}_{6}$  &  --- & $(24,60,30)$& $2$&$[hr,sl]$&$R4.2$&5 \\
\hline
$\{6,5\}_{4}$  &  --- & $(24,60,20)$ & $2$&$[hr,sr]$&$R9.16^*$&5\\
\hline
$\{6,5\}_{10}$  &  Petrie-dual of Small& $(24,60,20)$&$2$& $[hr,hl]$&$R9.15^*$ &6\\
&Stellated Dodecahedron &&&&& \\
\hline
$\{10,5\}_{6}$ &  Small Stellated & $(24,60,12)$ &$2$&$[hr,hr]$&$R13.8^*$&6 \\
&Dodecahedron&&&&&\\
\hline
$\{6,5\}_{10}$  &  Icosahedron & $(24,60,20)$&$1$&$[hr,hr]$ &$R9.15^*$&6 \\
\hline
$\{10,5\}_{6}$  &  Petrie-dual of & $(24,60,12)$& $1$&$[hr,hl]$&$R13.8^*$&6 \\
&Icosahedron&&&&&\\
\hline
$\{6,5\}_{10}$ &  Petrie-dual of  & $(24,60,20)$ &$1$& $[sr,sl]$&$R9.15^*$&7 \\
&Great Dodecahedron &&&&& \\
\hline
$\{10,5\}_{6}$  &  Great Dodecahedron& $(24,60,12)$ &$1$&$[sr,sr]$&$R13.8^*$&7\\
\hline
$\{6,5\}_{10}$  &  Great Icosahedron& $(24,60,20)$& $2$&$[sr,sr]$&$R9.15^*$&7 \\
\hline
$\{10,5\}_{6}$  &  Petrie-dual of & $(24,60,12)$&$2$&$[sr,sl]$&$R13.8^*$&7 \\
&Great Icosahedron&&&&&\\
\hline
\end{tabular}
\caption{The regular polyhedra of index $2$ with two vertex orbits.}
\label{tabone}}
\end{table}

A summary of these (families of) polyhedra is given in Table~\ref{tabone}, where the rows are in the same order as the diagrams in Figures~\ref{figpone} to \ref{figpsix}. The first column lists the type $\{p,q\}_t$ of a polyhedron, where $t$ is the length of the Petrie polygon, but note that this notation does not imply that $P$ is isomorphic to the universal regular map $\{p, q\}_t$. The last column refers to the complete enumeration, up to duality, of the orientable and non-orientable regular maps of small genus by Conder~\cite{con} (see also \cite{condob}), and gives the name of the underlying abstract polyhedron as a map in the listing.  In particular, $Rg.k$ designates the $k^{th}$ map among the orientable (\underbar{r}eflexibly) regular maps of genus $g$, in the labeling of \cite{con}, and $Rg.k^*$ denotes its dual. Thus $R4.2$ is the 2nd orientable regular map of genus $4$ in the listing. As the listing only begins with genus $2$, Table~\ref{tabone} has no entries for the spherical and toroidal maps $\{4,3\}_6$ and $\{6,3\}_4$, respectively, occurring in the first four rows of the last column. 

Moreover, the four families of polyhedra of types $\{4,5\}_{6}$ or $\{6,5\}_{4}$ are not generated from any ordinary regular polyhedron in the manner specified by Lemma~\ref{index12}, and so the second column has no entries at the corresponding positions. For edge length $1$, the two families each contain exactly one polyhedron with planar faces. These two polyhedra are the only planar-faced regular polyhedra of index $2$ with vertices on two orbits and already occur among the five polyhedra described by Wills~\cite{wills}. We shall see in \cite{cut2} that the two maps $\{4,5\}_{6}$ or $\{6,5\}_{4}$ for these four families of ``sporadic" polyhedra are dual to the maps of two regular polyhedra of index $2$ with vertices on just one orbit. 

It is interesting to note that Table~\ref{tabone} provides three separate instances of four regular, but geometrically dissimilar, polyhedra having the same regular map, namely those of types $\{10,3\}_{10}$, $\{6,5\}_{10}$, and $\{10,5\}_{6}$.

Recall that polyhedra in the same family are geometrically similar as they vary only in the value of $\lambda$, the ratio of the vertex orbit diameters, and we stipulate that $\lambda > 0$. If $\lambda$ is allowed to be negative, that is, if we reflect $S^\diamond$ about $o$, we get polyhedra that are still combinatorially equivalent, but are no longer geometrically similar. Such polyhedra belong to two different families, which nevertheless have the same underlying map because of the combinatorial equivalence.

The diagrams in Figure~\ref{figpone} and Figures~\ref{figpthree}, \ref{figpfour}, \ref{figpfive}, and \ref{figpsix} have been arranged so that structures in the same vertical column belong to families related in this way, and so share the same map, and structures in the same horizontal row are Petrie-duals of each other. Figure~\ref{figptwo} contains only two diagrams, since for octahedral symmetry the structures in the related families are also Petrie duals of each other. In this way the $22$ families of regular polyhedra of index 2 with vertices on two orbits fall naturally into groupings of four (or, in the case of octahedral symmetry, of two).

\end{document}